\documentclass[preprint,5p,twocolumn]{elsarticle}

\biboptions{sort&compress}

\makeatletter
\def\ps@pprintTitle{%
  \let\@oddhead\@empty
  \let\@evenhead\@empty
	\def\@oddfoot{\emph{Preprint}\hfil\emph{\today}}
  \let\@evenfoot\@oddfoot
}
\makeatother

\usepackage{amstext, amsmath, amssymb, graphicx, array, epsfig, psfrag}
\usepackage{color}
\usepackage{subfigure}
\usepackage{booktabs, multirow} 

\newcommand{\Vang}{\mathcal{V}^{\text{angle}}}
\newcommand{\Vmid}{\mathcal{V}^{\text{mid}}}
\newcommand{\V}{\mathcal{V}}
\newcommand{\bftcdotnabla}{\left(\bft \cdot \nabla\right)}
\newcommand{\bbR}{{\mathbb{R}}}

\newcommand{\bfM}{\boldsymbol M}
\newcommand{\bfT}{\boldsymbol T}
\newcommand{\bfn}{\boldsymbol n}
\newcommand{\bft}{\boldsymbol t}
\newcommand{\bfb}{\boldsymbol b}
\newcommand{\bfkappa}{\boldsymbol\kappa}

\newcommand{\bfp}{\boldsymbol p}

\newcommand{\bfr}{\boldsymbol r}

\newcommand{\bfN}{\boldsymbol N}
\newcommand{\bftheta}{\boldsymbol\theta}
\newcommand{\bfeta}{\boldsymbol\eta}
\newcommand{\bfv}{\boldsymbol v}
\newcommand{\bfgamma}{\boldsymbol\gamma}

\newcommand{\bff}{\boldsymbol f}
\newcommand{\bfu}{\boldsymbol u}
\newcommand{\bfeps}{\boldsymbol\varepsilon}
\newcommand{\bfsig}{\boldsymbol\sigma}
\newcommand{\bfzeta}{\boldsymbol\zeta}
\newcommand{\bfzero}{\boldsymbol 0}

\newcommand{\bfI}{\boldsymbol I}
\newcommand{\bfP}{\boldsymbol P}
\newcommand{\bfS}{\boldsymbol S}
\newcommand{\bfQ}{\boldsymbol Q}
\newcommand{\bfx}{\boldsymbol x}
\newcommand{\bfy}{\boldsymbol y}

\numberwithin{equation}{section}

\begin{document}

\begin{frontmatter}
\title{Variational Formulation of Curved Beams in Global Coordinates}
\author[jpk]{P. Hansbo}
\ead{peter.hansbo@jth.hj.se}

\author[umu]{M.G.\ Larson}
\ead{mats.larson@math.umu.se}

\author[umu]{K. Larsson\corref{cor1}}
\ead{karl.larsson@math.umu.se}

\cortext[cor1]{Corresponding author}

\address[jpk]{Department of Mechanical Engineering, J\"onk\"oping University, 
SE-55111 J\"onk\"oping, Sweden}
\address[umu]{Department of Mathematics and Mathematical Statistics, Ume{\aa} University, SE-90187 Ume{\aa}, Sweden}

\begin{abstract}
In this paper we derive a variational formulation for a linear curved beam which is natively expressed in global Cartesian coordinates. During derivation the beam midline is assumed to be implicitly described by a vector distance function which eliminates the need for local coordinates.
The only geometrical information appearing in the final expressions for the governing equations is the tangential direction, and thus there is no need to introduce normal directions along the curve.
As a consequence zero or discontinuous curvature, for example at inflection points, pose no difficulty in this formulation. Kinematic assumptions encompassing both Timoshenko and Euler--Bernoulli beam theories are considered.

With the exception of truly three dimensional formulations, models for curved beams found in literature are typically derived in the Frenet frame defined by the geometry of the beam midline.
While it is intuitive to formulate curved beam models in these local coordinates, the Frenet frame suffers from ambiguity and sudden changes of orientation in straight sections of the beam.

Based on the variational formulation we implement finite element models using global Cartesian degrees of freedom and discuss curvature coupling effects and locking. Numerical comparisons with classical solutions for both straight and curved cantilever beams under a tip load are given, as well as numerical examples illustrating curvature coupling effects.
\end{abstract}


\end{frontmatter}

\section{Introduction}

Models of one-dimensional elastic objects in $\bbR^3$, such as beams or rods with curved geometries,
are often established using a local equilibrium equation.
Such equilibrium equations are formulated using physical quantities, i.e. internal forces and moments,  defined in a local coordinate system, usually the Frenet frame which is defined through the differential geometry
of the curve by the Serret--Frenet formulas, cf., e.g., \cite{Reissner1962,Reissner1973,Simo1985,Reddy1993}. 

In this paper we instead establish a model for a curved beam using an equilibrium equation expressed in global Cartesian coordinates, and naturally the resulting governing equations are also expressed in global coordinates. The only geometrical information required for this formulation is the tangential direction, and thus zero curvatures do not pose any problems.
While a beam element defined in global coordinates was proposed in \cite{Gimena2008}, that formulation still depend on the Frenet frame to transform the differential equations to global coordinates.

Analogous descriptions have been formulated for two-dimensional elastic objects in $\bbR^3$, i.e. models for thin-shell structures.
Such models are often established using differential geometry to define the ruling differential equations \cite{Ci00}.
In the mid '90s, Delfour and Zol{\'e}sio \cite{DeZo95,DeZo96,DeZo97} instead established
elasticity models on surfaces using the signed distance function,
which can be used to describe the geometric properties of a surface.
In particular, the intrinsic tangential derivatives were
used for modeling purposes as the main differential geometric tool and the governing partial differential equations were
established in global Cartesian coordinates.
For one-dimensional objects a corresponding geometrical description can be formulated using a vector distance function,
which will give a formulation of the partial differential equations in three dimensions
based on the intrinsic differential operator which in the one-dimensional case is the tangential derivative.
While the term `intrinsic' is somewhat confusing in the case of one-dimensional curves we throughout this
paper continue to denote this approach intrinsic modeling.

The purpose of this paper is to begin to explore the possibilities of the intrinsic approach in finite element modeling of one-dimensional elastic structures embedded in $\mathbb{R}^3$, focusing on small strain Timoshenko and Euler--Bernoulli beam models.
Using this approach we are able to formulate a weak form of the governing partial differential
equations expressed in global coordinates, avoiding the need to introduce normal directions along the curve, by e.g. the Serret-Frenet formulas.
Rather, the only property inherent to a one-dimensional curve embedded in $\bbR^3$ which appears is the tangent direction.
In the derivation of the weak formulation of our beam problem we assume that the geometry
is defined using a vector distance function. However, in actual implementation of 
a finite element method 
any geometrical description may be used as long as it encompasses the tangent direction along the curve.



The remainder of this paper is dispositioned as follows.
In Section \ref{section:intrinsicgeneral} we present the intrinsic approach to codimension-two modeling,
i.e. using an implicit geometry description by a vector distance function. We relate this to the classical geometry description by a parametrized curve.
In Section \ref{section:intrinsicbeam} we start with the equilibrium equation of three-dimensional linear elasticity and a kinematic assumption
and derive a weak form of the governing equations for a
curved beam expressed in three dimensions.
In Section \ref{section:beamtheories} we explain how the Euler--Bernoulli
beam theory is encompassed in the kinematic assumption and
in Section \ref{section:locking} we give some remarks on the constraints
imposed on the approximation spaces when the thickness $t\rightarrow 0$.
Notes on our finite element implementations
and numerical results
in the form of a convergence study of a plane circular arc beam 
and in the form of numerical examples illustrating curvature effects
are presented in Section \ref{section:numerical}.
Finally, we give some concluding remarks in Section \ref{section:conclusion}.

\section{Intrinsic codimension-two modeling} \label{section:intrinsicgeneral}

\subsection{Basic notation}

Let $\Sigma$ be a smooth line embedded in $\bbR^3$, with tangent vector $\bft$.  We let 
$\bfp:\bbR^3 \rightarrow \Sigma$ be the closest point mapping, i.e. $\bfp (\bfx ) = \bfy$ 
where $\bfy \in \Sigma$ minimizes the Euclidian norm $|  \bfx - \bfy |_{\bbR^3}$.
The vector distance function, i.e. the vector between $\bfx$
and $\bfp(\bfx)$, we denote by $\bfzeta(\bfx) = \bfx - \bfp(\bfx)$.
The line $\Sigma$ is assumed to be the centerline of a beam with cross section $A$, which we for simplicity assume is constant. More precisely, $\Sigma$ passes through the centroid of $A$ and the beam 
occupies the volume
\begin{equation}
\Sigma \times A = \{ \bfx \in \bbR^3 \, : \, \bfx \in A\left(\bfp(\bfx)\right) \}
\end{equation}

The  linear projector
$\bfP_\Sigma = \bfP_\Sigma(\bfx)$, onto the tangent line of $\Sigma$ at $\bfx\in\Sigma$, 
is given by
\begin{equation}
\bfP_\Sigma = \bft\otimes\bft 
\end{equation}
where $\otimes$ denotes exterior product. We shall also need the projection 
\begin{equation}
\bfQ_\Sigma = \bfI - \bfP_\Sigma
\end{equation} 
onto the cross section plane orthogonal to $\bft$.
Let the gradient of a vector be defined by $\nabla\bfv=\nabla\otimes\bfv$.
We note that there is a neighborhood $\mathcal{N}(\Sigma ) \subset \bbR^3$ of $\Sigma$ such that 
$\bfp$ is an injective mapping in $\mathcal{N}(\Sigma)$ and $\Sigma\times A \subset \mathcal{N}(\Sigma)$.
Then for any function $v$ defined on $\Sigma$ we define the extension,
also denoted by $v$ to $\mathcal{N}(\Sigma) \subset \bbR^d$ by
\begin{align}
v ( \bfx ) = v(\bfp (\bfx ) )
\label{extension}
\end{align}

\subsection{Intrinsic curve geometry}

In this setting the geometry of the curve is implicitly defined through the vector distance function $\bfzeta(\bfx)$
defined above. On $\Sigma$ we have
\begin{align}
\bftcdotnabla \bfzeta &= \bfzero 
\label{tdell}
\\
\left( \bfn_i \cdot\nabla \right) \bfzeta &= \bfn_i \quad\text{for any $\bfn_i \bot \bft$}
\label{ndell}
\end{align}
As $\bftcdotnabla \bfzeta = \left( \nabla\bfzeta \right)^\text{T} \bft$
and $\left( \bfn_i \cdot\nabla \right) \bfzeta = \left( \nabla\bfzeta \right)^\text{T} \bfn_i$ the matrix
$\left( \nabla\bfzeta \right)^\text{T}$ clearly is a matrix with eigenvalues $\{0,1,1\}$ and $\bft$ is the eigenvector
associated with the zero eigenvalue. Thereby $\bft$ is uniquely determined down to its sign, but this will pose no problem in the intrinsic beam formulation.

As pointed out in the introduction, while the vector distance function is used in the derivation of the governing differential equations the geometrical description
used in the resulting method is arbitrary as long as the tangent direction is defined.
Typically, there are much more convenient ways than vector distance functions to describe the geometry in actual implementations.

\subsection{Differential curve qeometry}
In classical descriptions of curved beams (see e.g. \cite{Reissner1962,Reissner1973,Simo1985,Reddy1993}) the geometry is typically defined through an arc length parametrized curve $\bfr(\ell)$.
Recall that the unit tangent vector
is given by $\bft = {d \bfr}/{d \ell}$ and the curvature vector $\bfkappa$ of the curve is given by
\begin{align}
\bfkappa = \frac{d^2 \bfr}{d \ell^2} = \frac{d \bft}{d \ell}
= \bftcdotnabla\bft 
= \kappa \bfn
\end{align}
where $\bfn$ is the principal unit normal and $\kappa$ is the curvature given by $\kappa = | \bfkappa |$. The unit vector $\bfb = \bft\times\bfn$ is known as the binormal. We can also define the
curve torsion $\tau$ through
\begin{align}
\frac{d\bfb}{d\ell}=\bftcdotnabla\bfb=\bft\times\bftcdotnabla\bfn = -\tau \bfn
\end{align}
where we used that $\bftcdotnabla\bfb$ is
orthogonal to $\bfb$ and by the above cross product also is
orthogonal to $\bft$.
As $\bfn=\bfb\times\bft$ we may also write
\begin{align}
\frac{d\bfn}{d\ell}=\bftcdotnabla\bfn=
\bfb\times\bftcdotnabla\bft + \bftcdotnabla\bfb\times\bft
= -\kappa\bft + \tau\bfb
\end{align}

\begin{figure}[bt]
\begin{center}
\includegraphics{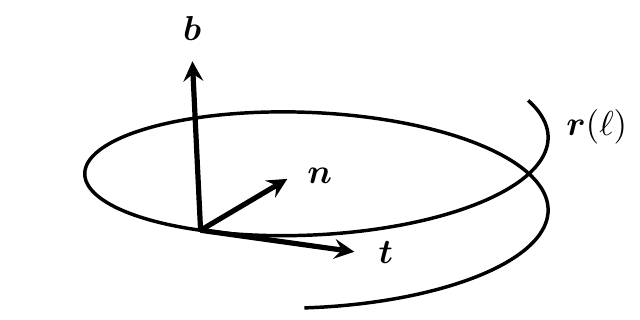}
\end{center}
\caption{Illustration of coordinate system given by the Frenet--Serret formulas}
\label{fig:serret}
\end{figure}
The above formulas define an orthonormal coordinate system $\{\bft,\bfn,\bfb\}$ along the curve (see Figure \ref{fig:serret}), and
 can be summarized in the following identity
known as the Frenet--Serret formulas
\begin{align}
\bftcdotnabla
\begin{bmatrix}
\bft \\
\bfn \\
\bfb
\end{bmatrix}
=
\begin{bmatrix}
0 & \kappa & 0 \\
-\kappa & 0 & \tau \\
0 & -\tau & 0
\end{bmatrix}
\begin{bmatrix}
\bft \\
\bfn \\
\bfb
\end{bmatrix}
\end{align}
Note that torsion is well defined only when $\kappa \neq 0$, so that $\bfn$ is defined.
For this reason special care have to be taken when using curves including zero curvatures
in formulations which depend on these normal directions to be defined.

\section{The intrinsic curved beam} \label{section:intrinsicbeam}

To explore the possibilities of the intrinsic approach for modeling one-dimensional
elastic structures embedded in $\mathbb{R}^3$ we in this section derive the
governing equations for a linear curved beam from the equations of linear elasticity
and a kinematic assumption.

\subsection{Kinematic assumption}

Based on the assumption that plane cross sections orthogonal to the midline remain plane
after deformation we assume that the displacement takes the following form 
\begin{equation} \label{displacementfield}
\bfu = \bfu_{\text{mid}} + \bftheta \times \bfzeta  
\end{equation}
where $\bfu_{\text{mid}}:\Sigma \rightarrow \mathbb{R}^3$ is the deformation of the midline 
and $\bftheta:\Sigma \rightarrow \mathbb{R}^3$ is an angle representing an infinitesimal rotation.
By the extension \eqref{extension} these functions are constant in the normal plane.
We may decompose $\bftheta$ as
\begin{align}
\bftheta = \bfQ_\Sigma \bftheta + \bfP_\Sigma \bftheta = \bfQ_\Sigma \bftheta + \bft \theta_t
\end{align}
where the first term $\bfQ_\Sigma \bftheta$ describes bending and shear 
and the second term $\bft \theta_t$ describes twist about the tangential axis.
Depending on how we define $\bftheta$ this kinematic assumption may encompass
both Timoshenko and Euler--Bernoulli beam theories (see Section \ref{section:beamtheories}).

\subsection{Strain} \label{section:strain}
We introduce the symmetric Cauchy strain tensor $\bfeps$ as
\begin{align}
\bfeps(\bfu) =\frac{1}{2}\left(\nabla\bfu + (\nabla\bfu)^{\rm T}\right)
\end{align}

Inserting the assumed displacement field \eqref{displacementfield} we note
that\\ $\bfQ_\Sigma\bfeps\left(\bfu_{\text{mid}} + \bftheta \times \bfzeta\right)\bfQ_\Sigma=\bfzero$
by the following three arguments.
Firstly, $\bfu_{\text{mid}}$ is constant in any normal direction which gives 
$\bfQ_\Sigma\nabla\bfu_{\text{mid}}\bfQ_\Sigma=\bfzero$.
Secondly, by \eqref{ndell} we have that $\left(\bfn_i \cdot \nabla\right)\left(\bfQ_\Sigma \bftheta \times \bfzeta \right)$ will be a vector in $\bft$-direction for any $\bfn_i \bot \bft$ and thus
$\bfQ_\Sigma\nabla\left(\left(\bfQ_\Sigma \bftheta\right) \times \bfzeta \right)\bfQ_\Sigma=\bfzero$.
Thirdly, by \eqref{tdell} and \eqref{ndell} we have that
$\bfQ_\Sigma \nabla(\theta_t \bft \times \bfzeta) \bfQ_\Sigma = -\bfQ_\Sigma \, \text{Skew}\left(\theta_t \bft \right) \bfQ_\Sigma$
which is a skew-symmetric matrix and thus gives no contribution to the strain.
In conclusion, using the assumed displacement field we may decompose the strain tensor into
\begin{align}
\bfeps(\bfu)
= \bfeps^P_\Sigma(\bfu) + \bfeps^S_\Sigma(\bfu)
\end{align}
where the in-line strain tensor $\bfeps^P_\Sigma(\bfu)$
and the shear strain tensor $\bfeps^S_\Sigma(\bfu)$ are
given by
\begin{align}
\bfeps^P_\Sigma(\bfu)
&= \bfP_\Sigma\bfeps(\bfu)\bfP_\Sigma
\\
\bfeps^{S}_\Sigma(\bfu)
&= \bfQ_\Sigma\bfeps(\bfu)\bfP_\Sigma + \bfP_\Sigma\bfeps(\bfu)\bfQ_\Sigma
\end{align}

Due to this decomposition of the strain tensor we have the
following relations for contractions between strain tensor components
\begin{align}
\bfeps_\Sigma^P ( \bfu ):\bfeps_\Sigma^S ( \bfv )
&= \bfeps_\Sigma^S ( \bfu ):\bfeps_\Sigma^P ( \bfv ) = 0 \label{cont1} \\
\bfeps_\Sigma^P ( \bfu ):\bfeps_\Sigma^P ( \bfv )
&= \left(\bfeps_\Sigma^P ( \bfu ) \cdot \bft\right) \cdot \left(\bfeps_\Sigma^P ( \bfv ) \cdot \bft\right) \\
\bfeps_\Sigma^S ( \bfu ):\bfeps_\Sigma^S ( \bfv )
&= 2 \left(\bfeps_\Sigma^S ( \bfu ) \cdot \bft\right) \cdot \left(\bfeps_\Sigma^S ( \bfv ) \cdot \bft\right) \label{cont3}
\end{align}
and we note that it suffices to calculate expressions for $\bfeps_\Sigma^P(\bfu) \cdot\bft$ and
$\bfeps_\Sigma^S(\bfu) \cdot\bft$ instead of the full tensors to compute these contractions.
Inserting our assumed displacement field \eqref{displacementfield} the expressions for the strains are given by
\begin{align}
\bfeps_\Sigma^P ( \bfu_{\text{mid}} ) \cdot \bft
&=
\bfP_\Sigma \bftcdotnabla \bfu_{\text{mid}}
\label{ePu}
\\
\bfeps_\Sigma^P (\bftheta \times \bfzeta ) \cdot \bft
&=
\bfP_\Sigma \left( \bftcdotnabla \bftheta \times \bfzeta \right)
\label{ePth}
\\
2 \bfeps_\Sigma^S ( \bfu_{\text{mid}} ) \cdot \bft
&=
\bfQ_\Sigma \bftcdotnabla \bfu_{\text{mid}}
\label{eSu}
\\
2 \bfeps_\Sigma^S (\bftheta \times \bfzeta ) \cdot \bft
&=
\bfQ_\Sigma \left( \bftcdotnabla \bftheta  \times \bfzeta \right)
-\bftheta\times\bft
\label{eSth}
\end{align}
Derivations of these strain expressions are supplied in Appendix~\ref{app:strain}.

\subsection{Governing equations}

The equilibrium equation of three dimensional linear elasticity reads
\begin{align}
- \nabla \cdot \bfsig &= \bff  \quad\text{in $\Sigma\times A$} \label{equilibrium} \\
\bfsig \cdot \bfn &= \bfzero \quad\text{on $\Sigma\times \partial A$}
\end{align}
where $\bfsig$ is the symmetric Cauchy stress tensor and $\bff$ is the body force density which
we for simplicity assume is constant over any cross-section $A$.
Note that we have left the boundary conditions on beam endpoints $\partial\Sigma\times A$ undefined.
We will return to these boundary conditions in Section~\ref{bcsection}.

Assuming that the stress components given by $\bfQ_\Sigma \bfsig \bfQ_\Sigma$ are zero
and that $\bft$ defines the material symmetry direction,
the constitutive relationship for an orthotropic\footnote{In case different shear moduli in
some orthogonal normal directions $\bfn_1$ and $\bfn_2$ is desired, $G \bfQ_\Sigma$ may be replaced by
$G_1 \bfn_1 \otimes \bfn_1 + G_2 \bfn_2 \otimes \bfn_2$.}
linear elastic material is given by
\begin{align} \label{constitutive}
\bfsig(\bfu) = E \bfeps^P_\Sigma(\bfu) + 2 G \bfeps^S_\Sigma(\bfu)
\end{align}
where $E$ is the elastic modulus and $G$ is the shear modulus \cite{ZiTa05}.
As a consequence we also have the relationships
\begin{align}
\bfP_\Sigma \bfsig(\bfu) \cdot \bft &= E \bfeps^P_\Sigma(\bfu) \cdot \bft
\\
\bfQ_\Sigma \bfsig(\bfu) \cdot \bft &= 2 G \bfeps^S_\Sigma(\bfu) \cdot \bft
\end{align}

\subsection{Function spaces}

As we have not yet defined boundary conditions on the beam ends $\partial\Sigma\times A$ we for now
only indicate that the essential boundary conditions are satisfied in the function spaces below
and refer the reader to Section \ref{bcsection} for the actual expressions.
On $\Sigma$ we introduce the function spaces
\begin{align}
\Vmid &= \left\{ \bfv_{\text{mid}} \in \left[H^1\left(\Sigma\right)\right]^3 : \text{$\bfv_{\text{mid}}$ with essential b.c.} \right\}
\\
\Vang &= \left\{ \bfeta \in \left[H^1\left(\Sigma\right)\right]^3 : \text{$\bfeta$ with essential b.c.}  \right\}
\end{align}
and using the vector distance function $\bfzeta$ we introduce the function space
\begin{align}
\V &= \left\{ \bfv=\bfv_{\text{mid}}+\bfeta\times\bfzeta : \ \bfv_{\text{mid}}\in\Vmid \, , \ \bfeta\in\Vang \right\}
\end{align}
We also introduce corresponding function spaces $\left\{\Vmid_0,\Vang_0,\V_0\right\}$ where the essential boundary conditions are homogenous. By the extension \eqref{extension} we
for $\bfv_{\text{mid}} \in \Vmid$ and $\bfeta \in \Vang$
have $\left\{\bfv_{\text{mid}},\bfeta\right\} \in \left[H^1\left(\Sigma\times A\right)\right]^3$.

\subsection{Weak formulation}
Multiplying \eqref{equilibrium} with a test function
$\bfv=\bfv_{\text{mid}}+\bfeta\times\bfzeta$ where $\bfv_{\text{mid}}\in\Vmid_0$ and $\bfeta\in\Vang_0$,
integrating over the domain $\Sigma \times A$, and applying Green's formula we end up with
\begin{align}
\begin{aligned}
\label{weak1}
\int_\Sigma &\int_A \bfsig(\bfu) : \bfeps(\bfv) \, dA d\Sigma
\\&=
\int_\Sigma \int_A \bff \cdot \bfv \, dA d\Sigma
+
\left[ \int_{A} (\bfsig(\bfu) \cdot \bft) \cdot \bfv \, dA \right]_{\partial \Sigma}
\end{aligned}
\end{align}
where the last term is handled by boundary conditions and is discussed in Section~\ref{bcsection}.

Using the constitutive relationship \eqref{constitutive} and the properties of in-line and shear strains in contraction \eqref{cont1}-\eqref{cont3} we have
\begin{align}
\int_\Sigma &\int_A \bfsig(\bfu) : \bfeps(\bfv) \, dA d\Sigma \nonumber
\\ &
\begin{aligned}
&=
\int_\Sigma \int_A
\left( \bfsig(\bfu)\cdot\bft \right)
\\&\qquad\qquad\quad \ \
\cdot
\left( \bfeps^P_\Sigma(\bfv)\cdot\bft + 2\bfeps^S_\Sigma(\bfv)\cdot\bft \right)
\, dA d\Sigma
\end{aligned}
\\ &
\begin{aligned}
&=
\int_\Sigma \int_A
\left( \bfP_\Sigma \bfsig(\bfu)\cdot\bft \right)\cdot
\left( \bfeps^P_\Sigma(\bfv)\cdot\bft\right)
\\&\qquad\qquad +
\left( \bfQ_\Sigma \bfsig(\bfu)\cdot\bft \right)\cdot
\left( 2\bfeps^S_\Sigma(\bfv)\cdot\bft\right)
\, dA d\Sigma
\end{aligned}
\end{align}
which when we instert $\bfv=\bfv_{\text{mid}} + \bfeta \times \bfzeta$ yields the four terms
\begin{align}
&\int_\Sigma \int_A \bfsig(\bfu) : \bfeps(\bfv) \, dA d\Sigma
\\&\ =
\int_\Sigma \int_A
\left( \bfP_\Sigma \bfsig(\bfu)\cdot\bft \right)\cdot
\left( \bfeps^P_\Sigma(\bfv_{\text{mid}})\cdot\bft\right)
\, dA d\Sigma \label{bilin1}
\\&\ \quad+
\int_\Sigma \int_A
\left( \bfP_\Sigma \bfsig(\bfu)\cdot\bft \right)\cdot
\left( \bfeps^P_\Sigma(\bfeta \times \bfzeta)\cdot\bft\right)
\, dA d\Sigma \label{bilin2}
\\&\, \quad+
\int_\Sigma \int_A
\left( \bfQ_\Sigma \bfsig(\bfu)\cdot\bft \right)\cdot
\left( 2\bfeps^S_\Sigma(\bfv_{\text{mid}})\cdot\bft\right)
\, dA d\Sigma \label{bilin3}
\\&\, \quad+
\int_\Sigma \int_A
\left( \bfQ_\Sigma \bfsig(\bfu)\cdot\bft \right)\cdot
\left( 2\bfeps^S_\Sigma(\bfeta \times \bfzeta)\cdot\bft\right)
\, dA d\Sigma \label{bilin4}
\end{align}
As $\bfv_{\text{mid}}$ due to the extension \ref{extension} is constant
over any cross-section $A$ we for the first term \eqref{bilin1} have
\begin{align}
\int_\Sigma \int_A
& \left( \bfP_\Sigma \bfsig(\bfu)\cdot\bft \right)\cdot
\left( \bfeps^P_\Sigma(\bfv_{\text{mid}})\cdot\bft\right)
\, dA d\Sigma
\\ &
\begin{aligned}
&=
\int_\Sigma
\left( \bfP_\Sigma \int_A \bfsig(\bfu)\cdot\bft  \, dA \right)
\\&\qquad\qquad\qquad
\cdot
\left( \bfP_\Sigma \bftcdotnabla \bfv_{\text{mid}} \right)
\, d\Sigma
\end{aligned}
\\ &=
\int_\Sigma
\bfN_\Sigma(\bfu) \cdot
\bftcdotnabla \bfv_{\text{mid}}
\end{align}
where we identify $\bfN_\Sigma= \bfP_\Sigma \int_{A} \bfsig\cdot\bft\, dA$ as the axial force
and use \eqref{ePu} in the second equality.
By the same arguments we for the third term \eqref{bilin3} have
\begin{align}
\int_\Sigma \int_A
&\left( \bfQ_\Sigma \bfsig(\bfu)\cdot\bft \right)\cdot
\left( 2\bfeps^S_\Sigma(\bfv_{\text{mid}})\cdot\bft\right)
\, dA d\Sigma
\\&
\begin{aligned}
&=
\int_\Sigma
\left( \bfQ_\Sigma \int_A \bfsig(\bfu)\cdot\bft  \, dA \right)
\\&\qquad\qquad\qquad
\cdot
\left( \bfQ_\Sigma \bftcdotnabla \bfv_{\text{mid}} \right)
\, d\Sigma
\end{aligned}
\\ &=
\int_\Sigma
\bfS_\Sigma(\bfu) \cdot
\bftcdotnabla \bfv_{\text{mid}}
\, d\Sigma
\end{align}
where we use \eqref{eSu} and $\bfS_\Sigma= \bfQ_\Sigma \int_{A} \bfsig\cdot\bft\, dA$ is the shear force.

Using \eqref{ePth} in the second term \eqref{bilin2} gives us
\begin{align}
\int_\Sigma &\int_A
\left( \bfP_\Sigma \bfsig(\bfu)\cdot\bft \right)\cdot
\left( \bfeps^P_\Sigma(\bfeta \times \bfzeta)\cdot\bft\right)
\, dA d\Sigma
\\ &=
\int_\Sigma \int_A
\left( \bfP_\Sigma \bfsig(\bfu)\cdot\bft \right)\cdot
\left( \bfQ_\Sigma \bftcdotnabla \bfeta \right) \times \bfzeta
\, dA d\Sigma
\\ &=
\int_\Sigma \int_A
\left(\bfzeta\times\left( \bfP_\Sigma \bfsig(\bfu)\cdot\bft \right)\right)\cdot
\left( \bfQ_\Sigma \bftcdotnabla \bfeta \right)
\, dA d\Sigma
\\ &
\begin{aligned}
&=
\int_\Sigma
\left( \bfQ_\Sigma \int_A \bfzeta\times\left( \bfsig(\bfu)\cdot\bft \right) \, dA \right)
\\&\qquad\qquad\qquad\qquad\qquad\quad
\cdot\left( \bfQ_\Sigma \bftcdotnabla \bfeta \right)
\, d\Sigma
\end{aligned}
\\ &=
\int_\Sigma
\bfM_\Sigma(\bfu)
\cdot \bftcdotnabla \bfeta
\, d\Sigma
\end{align}
where we identify $\bfM_\Sigma=\bfQ_\Sigma \int_A \bfzeta\times\left( \bfsig\cdot\bft \right) \, dA$ as the bending moment.
Similary, using \eqref{eSth} we for the fourth term \eqref{bilin4} have
\begin{align}
\int_\Sigma \int_A
&\left( \bfQ_\Sigma \bfsig(\bfu)\cdot\bft \right)\cdot
\left( 2\bfeps^S_\Sigma(\bfeta \times \bfzeta)\cdot\bft\right)
\, dA d\Sigma
\\ &
\begin{aligned}
&=
\int_\Sigma \int_A
\left( \bfQ_\Sigma \bfsig(\bfu)\cdot\bft \right)
\\&\qquad\qquad\qquad
\cdot
\left( \bfP_\Sigma \bftcdotnabla \bfeta \right) \times \bfzeta
\, dA d\Sigma
\\ &\quad
-\int_\Sigma \int_A
\left( \bfQ_\Sigma \bfsig(\bfu)\cdot\bft \right)\cdot
\bfeta\times\bft
\, dA d\Sigma
\end{aligned}
\\ &
\begin{aligned}
&=
\int_\Sigma 
\left( \bfP_\Sigma \int_A \bfzeta\times\left( \bfsig(\bfu)\cdot\bft \right) \, dA \right)
\\&\qquad\qquad\qquad\qquad\quad
\cdot
\left( \bfP_\Sigma \bftcdotnabla \bfeta \right)
\, d\Sigma
\\ &\quad
-\int_\Sigma
\left( \bfQ_\Sigma \int_A \left( \bfsig(\bfu)\cdot\bft \right) \, dA \right)
\cdot \bfeta\times\bft
\, d\Sigma
\end{aligned}
\\ &
=
\int_\Sigma \bfT_\Sigma(\bfu)\cdot
\bftcdotnabla \bfeta 
-\bfS_\Sigma(\bfu)\cdot
\bfeta \times \bft
\, d\Sigma
\end{align}
where $\bfT_\Sigma=\bfP_\Sigma \int_A \bfzeta\times\left( \bfsig\cdot\bft \right) \, dA$ is the torsion and recall that $\bfS_\Sigma$ is the shear force.

Collecting the above terms we have
\begin{align}
\int_\Sigma \int_A &\bfsig(\bfu) : \bfeps(\bfv) \, dA d\Sigma \nonumber
\\ &=
\int_\Sigma
\bfN_\Sigma(\bfu) \cdot
\bftcdotnabla \bfv_{\text{mid}}
\, d\Sigma
\\ &\quad+
\int_\Sigma
\bfS_\Sigma(\bfu) \cdot
\left(\bftcdotnabla \bfv_{\text{mid}} - \bfeta \times \bft \right)
\, d\Sigma
\\ &\quad+
\int_\Sigma
\bfM_\Sigma(\bfu) \cdot \bftcdotnabla \bfeta
\, d\Sigma
\\ &\quad+
\int_\Sigma
\bfT_\Sigma(\bfu) \cdot \bftcdotnabla \bfeta
\, d\Sigma
\end{align}
and note that the first term is associated with stretching, the second term is associated with shearing, the third term is associated with bending, and the last term is associated with twisting.
We will return to the actual expressions for the forces and moments in Section \ref{forcesandmoments}.

Before we present the problem in a more abstract setting with a bilinear form and linear functional
we in the next section turn to the boundary term of \eqref{weak1} and discuss suitable boundary conditions.

\subsection{Boundary conditions} \label{bcsection}
Consider the integral which appears in the boundary term of \eqref{weak1}.
Inserting $\bfv=\bfv_{\text{mid}}+\bfeta\times\bfzeta$ gives
\begin{align}
\int_{A} (\bfsig(\bfu) &\cdot \bft) \cdot \bfv \, dA
\nonumber\\&=
\int_{A} (\bfsig(\bfu) \cdot \bft) \cdot (\bfv_{\text{mid}}+\bfeta\times\bfzeta) \, dA
\\&
\begin{aligned}
&=
\left(\int_{A} (\bfsig(\bfu) \cdot \bft) \, dA \right) \cdot \bfv_{\text{mid}}
\\&\quad+
\left(\int_{A} \bfzeta\times(\bfsig(\bfu) \cdot \bft) \, dA \right) \cdot \bfeta
\end{aligned}
\end{align}
and by identifying forces and moments the boundary term for the beam ends can be written
\begin{align}
&\left[ \int_{A} (\bfsig(\bfu) \cdot \bft) \cdot \bfv \, dA \right]_{\partial \Sigma} \nonumber
\\
&\quad=
\bigl[
\left( \bfN_\Sigma(\bfu) + \bfS_\Sigma(\bfu) \right) \cdot \bfv_{\text{mid}} \nonumber
\\
&
\qquad\qquad\qquad\qquad
+  \left(\bfM_\Sigma(\bfu) + \bfT_\Sigma(\bfu) \right) \cdot\bfeta
\bigr]_{\partial \Sigma}
\\
&\begin{aligned}
&\quad=
\bigl[
\bfN_\Sigma(\bfu)\cdot\bft v_t
+ \bfS_\Sigma(\bfu) \cdot(\bfQ_\Sigma \bfv_{\text{mid}})
\\
&\qquad\qquad\quad
+ \bfM_\Sigma(\bfu) \cdot(\bfQ_\Sigma \bfeta)
+ \bfT_\Sigma(\bfu) \cdot\bft \eta_t
\bigr]_{\partial \Sigma}
\end{aligned}
\end{align}
We readily see that the natural and essential boundary conditions are those listed in
Table~\ref{bctable}. Note that we have to select one boundary condition for each row
of the table. 
We will not list all possible combinations of boundary conditions but remark that
if we select only essential boundary conditions we have a fixed (clamped) boundary and
if we select only natural boundary conditions we have a free boundary.

\begin{table}[ht]
\centering
\begin{tabular}{ l  c  c }
\toprule
    & Natural condition & Essential condition \\
	\midrule                        
  Streching & $\bfN_\Sigma\cdot\bft=\overline{N}$ & $u_t=\overline{u}_t$ \\
	\midrule
  Shearing & $\bfS_\Sigma=\overline{\bfS}$ & $\bfQ_\Sigma \bfu_{\text{mid}}=\overline{\bfu}_\bot$ \\
	\midrule
	Bending & $\bfM_\Sigma=\overline{\bfM}$ & $\bfQ_\Sigma \bftheta=\overline{\bftheta}_\bot$ \\
	\midrule
  Twisting & $\bfT_\Sigma\cdot\bft=\overline{T}$ & $\theta_t=\overline{\theta}_t$ \\
  \bottomrule
\end{tabular}
\caption{Associated natural and essential boundary conditions for stretching, shearing, bending, and twisting}
\label{bctable}
\end{table}

\subsection{Forces and moments} \label{forcesandmoments}
In this section we collect expressions for
the axial force $\bfN_\Sigma:= \bfP_\Sigma \int_{A} \bfsig(\bfu)\cdot\bft\, dA$, transverse shear force $\bfS_\Sigma:= \bfQ_\Sigma\int_{A} \bfsig(\bfu)\cdot\bft\, dA$, 
bending moment $\bfM_\Sigma:= \bfQ_\Sigma \int_{A}\bfzeta\times \left(\bfsig(\bfu)\cdot\bft\right) \, dA$, and torsion $\bfT_\Sigma:= \bfP_\Sigma \int_{A}\bfzeta\times \left(\bfsig(\bfu)\cdot\bft\right) \, dA$.
Introducing the area $|A|$, the tensor of area moments of inertia $\bfI_\Sigma$, and the polar inertia $J_\Sigma$;
\begin{align}
|A| &:= \int_A dA \label{area}
\\
\bfI_\Sigma &:= \int_A \bfzeta \times \bft \otimes \bfzeta \times \bft \, dA \label{areainert}
\\
J_\Sigma &:= \int_A \bfzeta \cdot \bfzeta \, dA \label{polar}
\end{align}
we have the following expressions for the forces and moments
\begin{align}
\label{normalforce}
\bfN_\Sigma(\bfu_{\text{mid}}) &=
E |A| \bfP_\Sigma \bftcdotnabla \bfu_{\text{mid}}
\\
\label{shearforce}
\bfS_\Sigma(\bfu_{\text{mid}}, \bftheta) &=
G |A| \left( \bfQ_\Sigma \bftcdotnabla \bfu_{\text{mid}} - \bftheta \times\bft \right)
\\
\label{bendingmoment}
\bfM_\Sigma(\bftheta) &=
E \bfI_\Sigma \bftcdotnabla \bftheta
\\
\label{torsion}
\bfT_\Sigma(\bftheta) &= 
G J_\Sigma \bfP_\Sigma \bftcdotnabla \bftheta
\end{align}
and we supply the derivations of these expressions in Appendix \ref{forcesandmomentsderivation}.
In Section \ref{section:curvature} we give alternate expressions for these forces and moments which
more explicitly state their dependence on the curvature.
Note that we have not included correction factors in the expressions for shear force and torsion
to account for the difference between the assumed displacement field and the actual deformation
of the cross-section.

\subsection{Abstract weak form} \label{abstractweakform}
We now collect the results of the above sections
in the following
abstract form of the weak problem:
Find $\bfu \in \V$ such that
\begin{align} \label{bilinear}
a(\bfu;\bfv) = l(\bfv) \qquad \text{for all $\bfv \in \V_0$}
\end{align}
where the bilinear form $a(\cdot;\cdot)$ is given by
\begin{align}
\begin{aligned}
a(\bfu;\bfv)
&=
a_{\text{Stretch}}( \bfu_{\text{mid}} ; \bfv_{\text{mid}} )
+
a_{\text{Shear}}(\bftheta;\bfeta)
\\&\quad
+
a_{\text{Bend}}( \bfu_{\text{mid}}, \bftheta ; \bfv_{\text{mid}}, \bfeta )
+
a_{\text{Twist}}(\bftheta;\bfeta)
\end{aligned}
\end{align}
where 
\begin{align}
a_{\text{Stretch}}
&=
\int_\Sigma
E |A| \left( \bfP_\Sigma \bftcdotnabla \bfu_{\text{mid}}) \right) \nonumber
\\&\qquad\qquad\qquad\quad
\cdot
\left( \bfP_\Sigma \bftcdotnabla \bfv_{\text{mid}})  \right)
\, d\Sigma
\label{aStrech}
\\ 
a_{\text{Shear}}
&=
\int_\Sigma
G |A| \left( \bfQ_\Sigma \bftcdotnabla \bfu_{\text{mid}} - \bftheta \times\bft \right)
\nonumber \\ &\qquad\qquad\quad
\cdot
\left(\bfQ_\Sigma \bftcdotnabla \bfv_{\text{mid}} - \bfeta \times \bft \right)
\, d\Sigma
\label{aShear}
\\
a_{\text{Bend}}
&=
\int_\Sigma
\left( E \bfI_\Sigma  \bftcdotnabla \bftheta \right)
\cdot
\left( \bftcdotnabla \bfeta \right)
\, d\Sigma
\label{aBend}
\\
a_{\text{Twist}}
&=
\int_\Sigma
G J_\Sigma \left( \bfP_\Sigma \bftcdotnabla \bftheta \right) \cdot
\left(\bfP_\Sigma \bftcdotnabla \bfeta \right)
\, d\Sigma
\label{aTwist}
\end{align}
As $\bff$ is assumed constant over any cross-section $A$ the linear functional $l(\cdot)$ is given by
\begin{align}
\begin{aligned}
l(\bfv) &= |A| \int_\Sigma \bff \cdot \bfv_{\text{mid}} \, d\Sigma
\\&\quad+\bigl[
\bfN_\Sigma(\bfu)\cdot\bft v_t
+ \bfS_\Sigma(\bfu) \cdot(\bfQ_\Sigma \bfv_{\text{mid}}) 
\\&\qquad\qquad
+ \bfM_\Sigma(\bfu) \cdot(\bfQ_\Sigma \bfeta)
+ \bfT_\Sigma(\bfu) \cdot\bft \eta_t
\bigr]_{\partial \Sigma}
\end{aligned}
\end{align}
with applicable natural boundary conditions (see Table~\ref{bctable}). For essential boundary conditions imposed on
$\bfu$ the corresponding homogenous boundary conditions are imposed on $\bfv$,
rendering any boundary term above without a natural boundary condition to be zero.

\section{Euler--Bernoulli beam theory} \label{section:beamtheories}
In Timoshenko beam theory the beam cross-section is assumed plane after deformation
but is allowed to shear from the beam midline. We may thus choose the approximation
of the beam midline $\bfu_{\text{mid}}$ independently from the angle $\bftheta$.
Apart from considerations needed to avoid locking effects when the beam is thin, this means that we
for a Timoshenko beam directly use the bilinear form given in the previous section.

Introducing the shear angle $\bfQ_\Sigma \bfgamma$ defined such that
\begin{align} \label{shearingangle}
\left(\bfQ_\Sigma \bfgamma \right) \times \bft = \bfQ_\Sigma \bftcdotnabla \bfu_{\text{mid}} -(\bfQ_\Sigma \bftheta) \times \bft
\end{align}
we may express the shear force in terms of the shear angle by
$\bfS_\Sigma =
G |A| \left(\bfQ_\Sigma \bfgamma \right)\times\bft$.

In Euler--Bernoulli beam theory the beam cross-section is assumed plane and orthogonal
to the beam midline after deformation. Letting the shearing angle $\bfQ_\Sigma \bfgamma \rightarrow \bfzero$ in \eqref{shearingangle} we have that
\begin{align}
(\bfQ_\Sigma \bftheta) \times \bft = \bfQ_\Sigma \bftcdotnabla \bfu_{\text{mid}}
\end{align}
which means that a suitable definition for
the angle $\bftheta$ is
\begin{equation}
\bftheta = \nabla \times \bfu_{\text{mid}} + \bft \theta_t
\end{equation}
This definition for $\bftheta$ in combination with \eqref{displacementfield} constitutes the
Euler--Bernoulli kinematic assumption.
The term $(\bfQ_\Sigma \bftheta) \times \bft$ can then be rewritten
\begin{align}
(\bfQ_\Sigma \bftheta) \times \bft
&=
-\bfQ_\Sigma\left( \bft\times(\bfQ_\Sigma \bftheta) \right)
\\&=
-\bfQ_\Sigma \left( \bft \times \left( \nabla \times \bfu_{\text{mid}} \right) \right)
\\ &=
-\bfQ_\Sigma \nabla \bfu_{\text{mid}} \cdot \bft + \bfQ_\Sigma \bftcdotnabla \bfu_{\text{mid}}
\\ &=
\bfQ_\Sigma \bftcdotnabla \bfu_{\text{mid}}
\end{align}
just as required when the shearing angle is zero.
By the Euler--Bernoulli kinematic assumption we thus have a shear force \eqref{shearforce}
which is zero.
Further, using this kinematic assumption we have that
\begin{align}
\bfQ_\Sigma \bftcdotnabla \bftheta
&=
\bfQ_\Sigma \bftcdotnabla \left(\bfQ_\Sigma \bftheta + \theta_t \bft \right)
\\&
\begin{aligned}
=
\bfQ_\Sigma ( &\bftcdotnabla \left(\bfQ_\Sigma \bftheta\right)
\\ &+ \theta_t \bftcdotnabla\bft + \bft \bftcdotnabla \theta_t )
\end{aligned}
\\&=
\bfQ_\Sigma \bftcdotnabla \left(\nabla \times \bfu_{\text{mid}}\right)
+ \theta_t \bfkappa
\end{align}
and
\begin{align}
\bfP_\Sigma &\bftcdotnabla \bftheta \nonumber \\
&=
\bfP_\Sigma \bftcdotnabla \left( \theta_t \bft + \bfQ_\Sigma \bftheta \right)
\\&
\begin{aligned}
&=
\bfP_\Sigma \left( \bft \bftcdotnabla \theta_t + \theta_t \bfkappa \right)
\\&\quad
+
\bfP_\Sigma
\left( \bftcdotnabla \left(\bfQ_\Sigma\right) \bftheta
+ \bfQ_\Sigma \bftcdotnabla\bftheta \right)
\end{aligned}
\\&=
\bft \left( \bftcdotnabla \theta_t
- \left(\nabla \times \bfu_{\text{mid}}\right) \cdot \bfkappa \right)
\end{align}
where we note that the first term contains second order derivatives of $\bfu_{\text{mid}}$
which implies that $\bfu_{\text{mid}} \in \left[H^2\left(\Sigma\right)\right]^3$ is needed for this term
to be well defined.

In conclusion, using the Euler--Bernoulli assumption that cross-sections normal to the midline
remain normal after deformation, \eqref{aShear} vanishes and \eqref{aBend} reduces to
\begin{align}
&a_{\text{Bend}}( \bfu_{\text{mid}},\theta_t ; \bfv_{\text{mid}},\eta_t )
\nonumber \\&\qquad
\begin{aligned}
&=
\int_\Sigma
\left( E \bfI_\Sigma \left( \bftcdotnabla \left(\nabla \times \bfu_{\text{mid}} \right) + \theta_t \bfkappa \right) \right)
\\ &\qquad\qquad\quad
\cdot\left( \bftcdotnabla \left(\nabla \times \bfv_{\text{mid}}\right) + \eta_t \bfkappa \right)
\, d\Sigma
\end{aligned}
\end{align}
and \eqref{aTwist} becomes
\begin{align}
&a_{\text{Twist}}( \bfu_{\text{mid}},\theta_t ; \bfv_{\text{mid}},\eta_t )
\nonumber \\&\qquad
\begin{aligned}
&=
\int_\Sigma
G J_\Sigma \left( \bftcdotnabla\theta_t - \left(\nabla \times \bfu_{\text{mid}}\right)\cdot\bfkappa \right)
\\ &\qquad\qquad\quad
\left(\bftcdotnabla\eta_t - \left(\nabla \times \bfv_{\text{mid}}\right)\cdot\bfkappa \right)
\, d\Sigma
\end{aligned}
\end{align}
The appropriate function space for the midline in the Euler--Bernoulli formulation is $\Vmid \cap \left[H^2\left(\Sigma\right)\right]^3$.
Note that the terms above in contrast to the formulation in Section~\ref{abstractweakform} explicitly contain
the curvature vector $\bfkappa$ which means
that the geometry needs to be defined in a manner such that the curvature is readily available. This is
typically not a problem in most implementations but it should however be noted that when we
include torsional effects in the Euler--Bernoulli beam theory, we lose one of the
benefits of the intrinsic beam formulation, i.e., that only information of the tangent is needed.
Still, there is no need for introducing normals along the curve and zero curvatures pose
no problem.

\section{Cartesian approximation spaces and locking effects} \label{section:locking}
In the formulation above it is natural to use approximation spaces where the degrees of freedom
are expressed in Cartesian coordinates.
The choice of approximation space is however nontrivial as there are compatibility
requirements introduced both by the curvature and by the mixed formulation in itself.

\subsection{Thickness scaling}
Assuming that the thickness of the beam is controlled by a parameter $t$
proportional to the diameter of the cross-section we have the following scaling
of the area, area moments of inertia, and polar inertia
\begin{align}
|A| = t^2 |A|^{\text{ref}}  \qquad
\bfI_\Sigma = t^4 \bfI_\Sigma^{\text{ref}}  \qquad
J_\Sigma = t^4 J_\Sigma^{\text{ref}}
\end{align}
Dividing all terms in $a(\bfu,\bfu)$ by $t^4$ and let $t\rightarrow 0$ we get constraints
\begin{align} \label{compreq}
\left(\bfQ_\Sigma \bftcdotnabla \bfu_{\text{mid}} - \bftheta \times \bft \right) &\rightarrow 0
\\ \label{compreq2}
\bfP_\Sigma \bftcdotnabla \bfu_{\text{mid}} &\rightarrow 0
\end{align}
While the first expression could be handled by using Euler--Bernoulli theory as suggested in
Section~\ref{section:beamtheories}.
We could also view the Euler--Bernoulli theory as a consequence of this expression
when $t\rightarrow 0$ and thus need to choose approximation space for $\bfu_{\text{mid}}$ such that
$\Vmid \in \left[H^2\left(\Sigma\right)\right]^3$ and the space $\Vang$
such that $\bftheta$ may fulfill this compatibility requirement.

However, the second expression still is problematic due
to the non-linear nature of $\bfP_\Sigma$. We consider this problem in the next section.

\subsection{Curvature locking and reduced integration}
While the expressions in \eqref{aStrech}-\eqref{aTwist} might look simple, the geometric
information contained in $\bft$, $\bfP_\Sigma$, and $\bfQ_\Sigma$ may be highly non-linear. Following the separation of tangential and cross-sectional components presented in Appendix~\ref{forcesandmomentsseparation} we get the following expressions
\begin{align}
\bfQ_\Sigma \bftcdotnabla \bfu_{\text{mid}}
&=
u_t \bfkappa 
+ \bfQ_\Sigma \bftcdotnabla \left(\bfQ_\Sigma\bfu_{\text{mid}}\right)
\\
\bft\cdot \bftcdotnabla \bfu_{\text{mid}}
&=
\bftcdotnabla u_t
- \left(\bfkappa\cdot\bfu_{\text{mid}}\right)
\\
\bfQ_\Sigma \bftcdotnabla \bftheta
&=
\theta_t \bfkappa 
+ \bfQ_\Sigma \bftcdotnabla \left(\bfQ_\Sigma\bftheta\right)
\\
\bft\cdot \bftcdotnabla \bftheta
&=
\bftcdotnabla \theta_t
-\left(\bfkappa\cdot\bftheta\right)
\end{align}
where we in each term note a mixture of derivatives and curvature
which introduce compatibility requirements, especially in the limit $t\rightarrow 0$
as noted above. As we interpolate $\bfu_{\text{mid}}$ and $\bftheta$ in Cartesian
coordinates it is not trivial to manipulate the interpolation spaces such that
the compatibility requirements are fulfilled. To avoid locking phenomenons
stemming from these incompatibilities we therefore resort to reduced quadrature
of terms \eqref{aStrech} and \eqref{aShear}.

\section{Numerical examples} \label{section:numerical}

The numerical results presented in this paper consist of a convergence study
where we compare tip deflection in two model problems with analytical results from
\cite{timoshenko1969theory}.
Further, to illustrate how the curvature couples
stretching, shearing, bending and twisting 
we include some numerical examples of deformation of
initially curved beams.

\subsection{Implementation}

We interpolate the vector valued functions $\bfu_{\text{mid}}$ and $\bftheta$ in Cartesian coordinates using piecewise polynomial interpolation with respect to the arc length $\ell$. The directional derivatives $\bftcdotnabla$ in this setting are easily evaluated. As indicated by the above formulas for the Timoshenko beam
it is sufficient to know the tangent vector $\bft$ at each quadrature point.
For convenience we choose finite element spaces which, in the straight case, do not exhibit locking when $t\rightarrow 0$.
Clearly, choosing  $C^1$-continuous piecewise cubic interpolation for $\bfu_{\text{mid}}$ and
continuous piecewise quadratic interpolation for $\bftheta$ defines compatible approximation spaces
in the straight case.
While less obvious, it is well known that choosing continuous piecewise quadratic interpolation
for $\bfu_{\text{mid}}$ and continuous piecewise linear interpolation for $\bftheta$ also
constitutes compatible approximation spaces.
This is due to the fact that there is a large enough $C^1(\Sigma)$ subspace within
the space of continuous piecewise quadratic functions.
Numerical results are provided for the following methods:
\begin{itemize}
\item Timoshenko beam using continuous piecewise quadratic interpolation for the midline and continuous piecewise linear interpolation for the angle (P2-P1).
\item Timoshenko beam with $C^1$-continuous cubic Hermite interpolation for the midline and continuous piecewise linear interpolation for the angle (H3-P2).
\end{itemize}

In our implementations the geometry was both represented exactly and by cubic Hermite spline interpolation.
The difference in the results was however small.

Integration was performed using Gauss quadrature points distributed along each element according
to the arc length. In the curved examples reduced integration of terms \eqref{aStrech} and \eqref{aShear} was used and
good results achieved when
the number of quadrature points were chosen such that polynomials up to an order of three was integrated exactly.

\subsection{Convergence study}
\begin{figure*}[tb]
\centering
\subfigure[Straight beam]{
\includegraphics[width=0.6\linewidth]{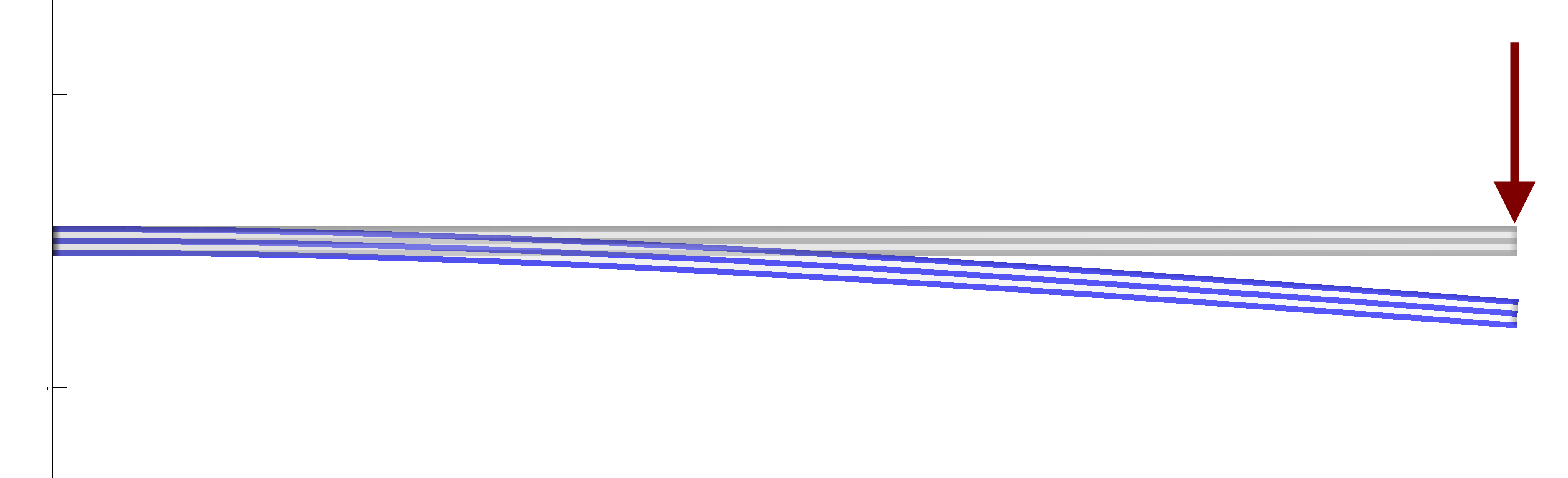}
\label{fig:straightbeam}
}
\subfigure[Curved beam]{
\includegraphics[width=0.3\linewidth]{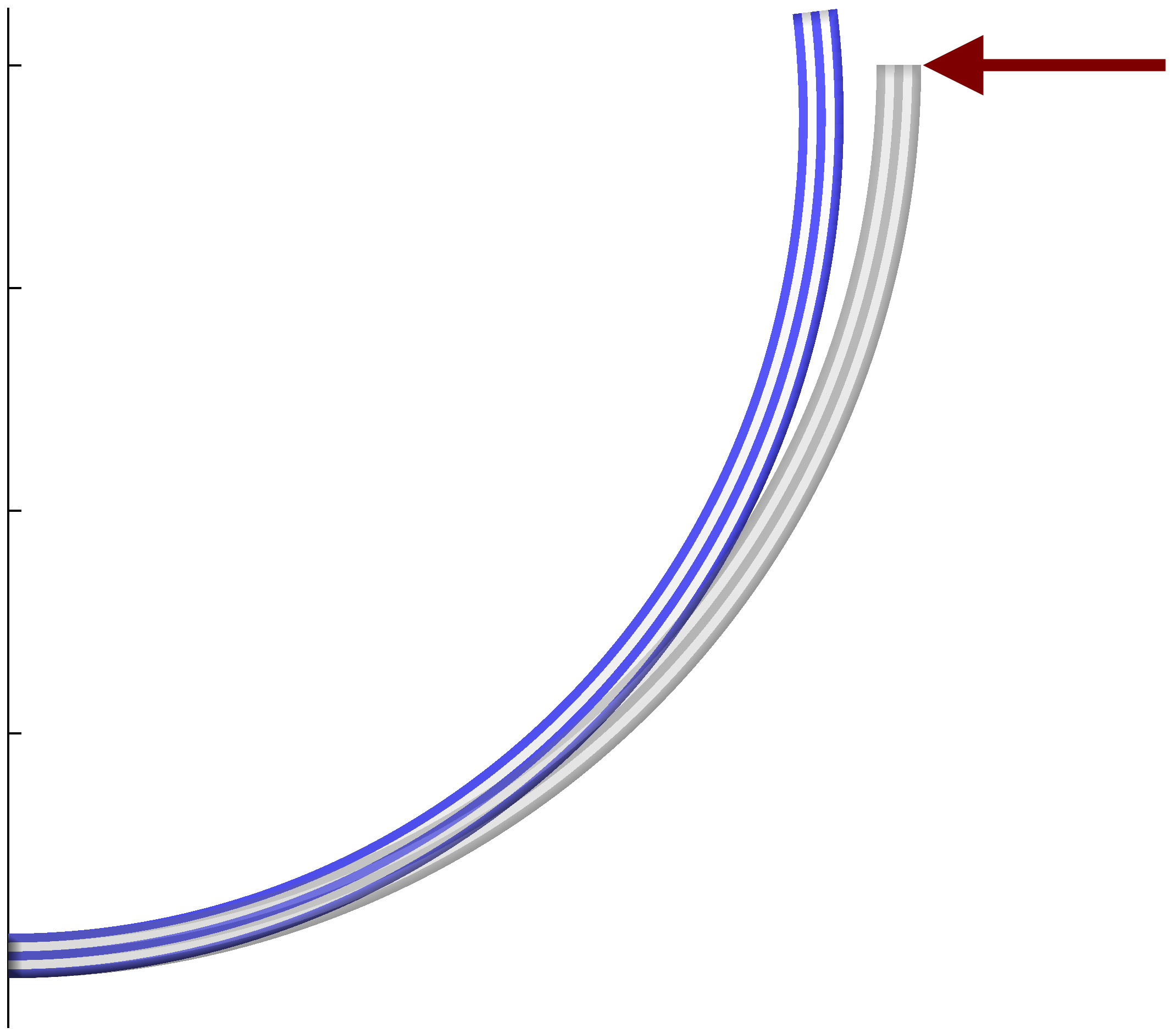}
\label{fig:curvedbeam}
}
\caption{Deflection of a straight (a) and a circular arc (b) cantilever beams under a tip load}
\label{fig:analytical}
\end{figure*}

In this paper we limit our convergence study to two
basic model problems presented in \cite{timoshenko1969theory}.
The first model problem is that of a straight cantilever beam of length $L$ and thickness $t$ under a tip load $P$ as illustrated in Figure~\ref{fig:straightbeam}. An analytical solution for the tip
deflection in $y$-direction for a beam of unit depth is given by
\begin{align}
u_y = -\frac{P}{6 E I} \left( (4+5\nu) \frac{t^2 L}{4} + (2 L^3) \right)
\end{align}
where $\nu$ is the Poisson's ratio. Assuming a rectangular cross-section the second moment of area is $I = t^3/12$.

The second model problem is that of a curved cantilever beam in the shape of a quarter circle under a tip load $P$ as illustrated in Figure~\ref{fig:curvedbeam}.
With inner radius $a$ and outer radius $b$ such that $b-a=t$
an analytical solution is given by
\begin{align}
u_x = -\frac{P \pi (a^2+b^2)}{E\left[ (a^2-b^2) + (a^2+b^2) \log(b/a) \right]}
\end{align}
where we again assume a rectangular cross-section and unit depth.

However, these analytical solutions are only exact if the forces at the end of the bar follow the same parabolic distributions as the shearing stress, a property which the implemented methods do not encompass. 
As noted in \cite{timoshenko1969theory, dur5010}, if this condition is not fulfilled the solutions
above are only approximations.
Hence we are comparing our numerical solutions with an analytical approximation of the solution which means that the numerical methods will eventually appear to stop converging in our convergence plots.

\begin{figure}
\centering
\subfigure[t=0.1]{
\includegraphics[width=1.0\linewidth]{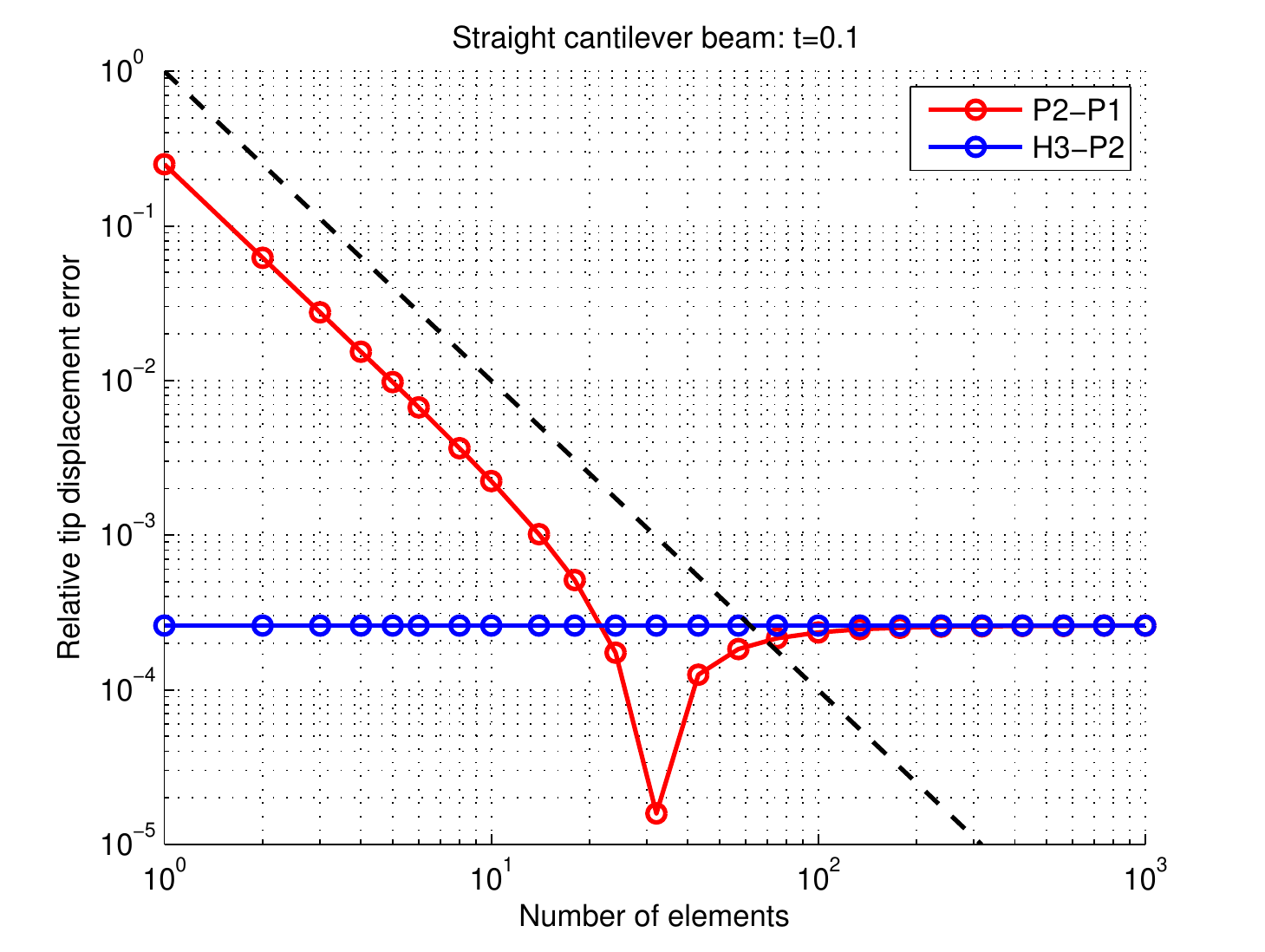}
\label{fig:straightconv-1}
}
\subfigure[t=0.001]{
\includegraphics[width=1.0\linewidth]{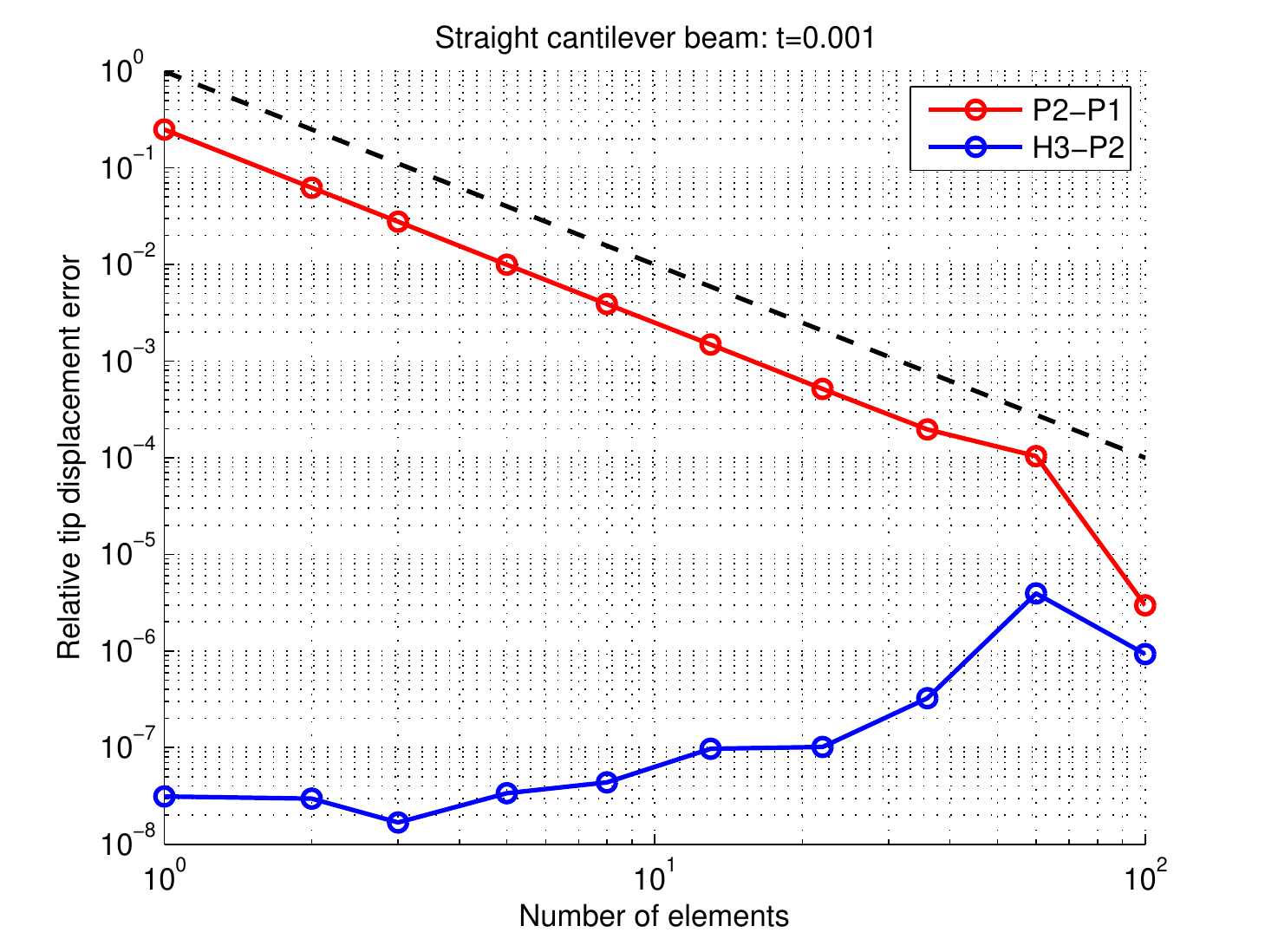}
\label{fig:straightconv-3}
}
\caption{Convergence of tip deflection for a straight cantilever beam under tip load which a dashed reference line corresponding to a convergence order of two}
\label{fig:straightconv}
\end{figure}

In the case of a straight beam expressions \eqref{aStrech}-\eqref{aTwist} are identical
with expressions for the irreducible weak form for the Timoshenko beam presented in \cite{ZiTa05}.
It is therefor unsurprising that the P2-P1 method gives the convergence behavior
depicted in Figure~\ref{fig:straightconv} which indicate a convergence order of 2.
Already with only one element, the H3-P2 method gives an approximation which is better than
the analytical approximation. Due to the choice of approximation spaces, neither method exhibits locking when $t\rightarrow 0$ as indicated by Figure \ref{fig:straightconv-3}.

\begin{figure}
\centering
\subfigure[t=0.1]{
\includegraphics[width=1.0\linewidth]{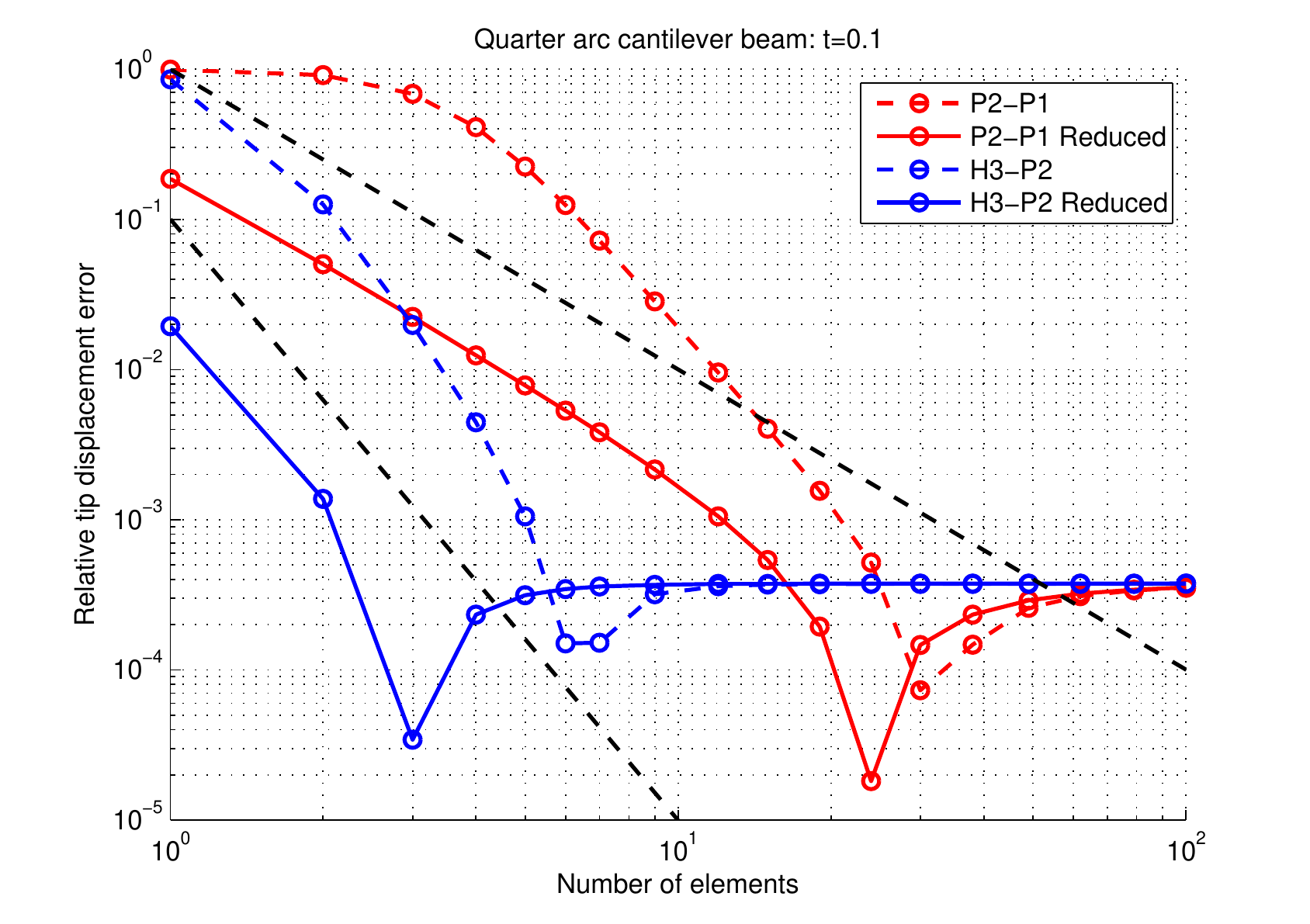}
\label{fig:curvedconv-1}
}
\subfigure[t=0.001]{
\includegraphics[width=1.0\linewidth]{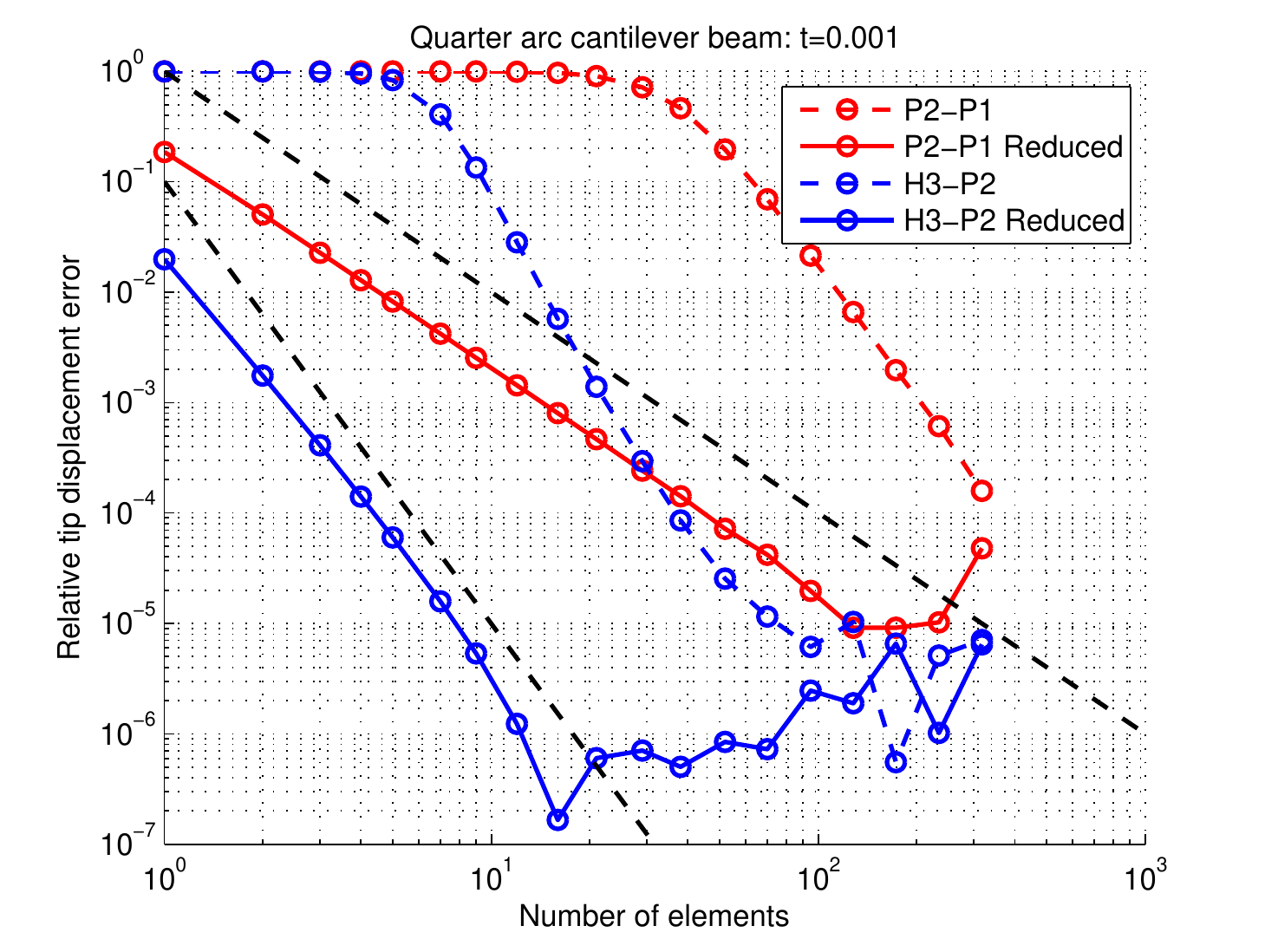}
\label{fig:curvedconv-3}
}
\caption{Convergence of tip deflection for a circular arc cantilever beam under tip load with dashed reference lines in black corresponding to a convergence order of two and four}
\label{fig:curvedconv}
\end{figure}

In the case of the quarter arc cantilever beam both methods are prone to locking as may be seen
in Figure \ref{fig:curvedconv} and clearly this effect increase with smaller $t$.
The use of reduced quadrature on terms \eqref{aStrech} and \eqref{aShear} however appears to successfully
remove the locking effects giving the P2-P1 method a convergence order of 2 and the H3-P2 method
a convergence order of 4.

\subsection{Curvature effects} \label{section:curvature}

\begin{figure}
\centering
\subfigure[Twisting moment]{
\includegraphics[height=0.4\linewidth]{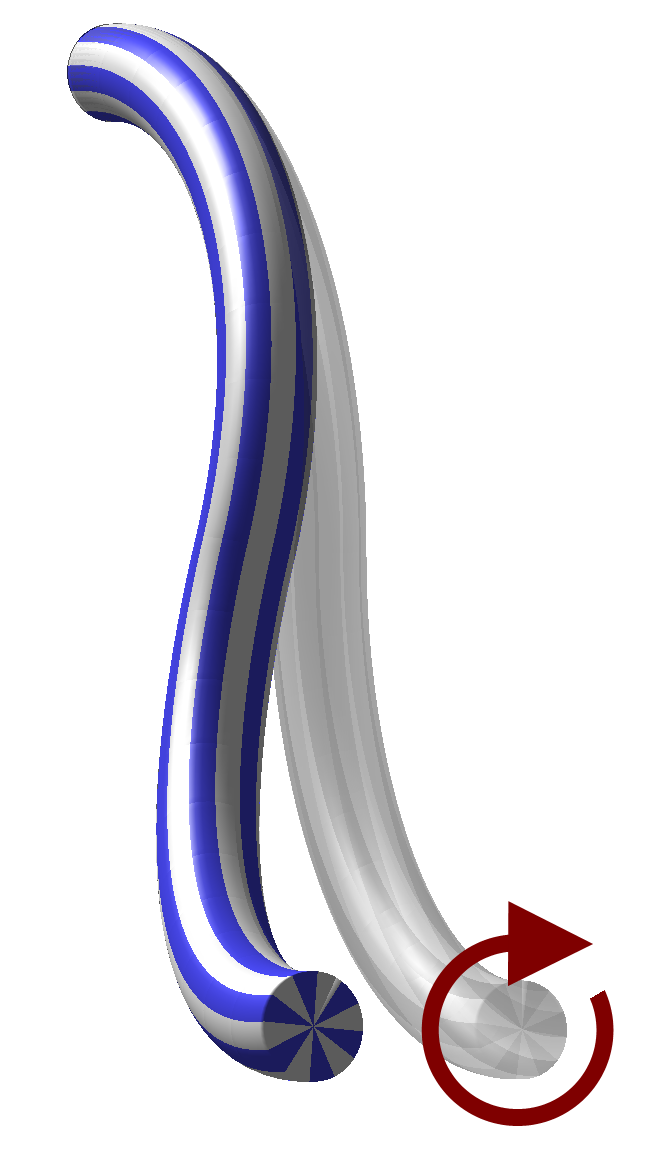}
\label{fig:sshape}
}
\subfigure[Point load orthogonal to curvature]{
\includegraphics[height=0.5\linewidth]{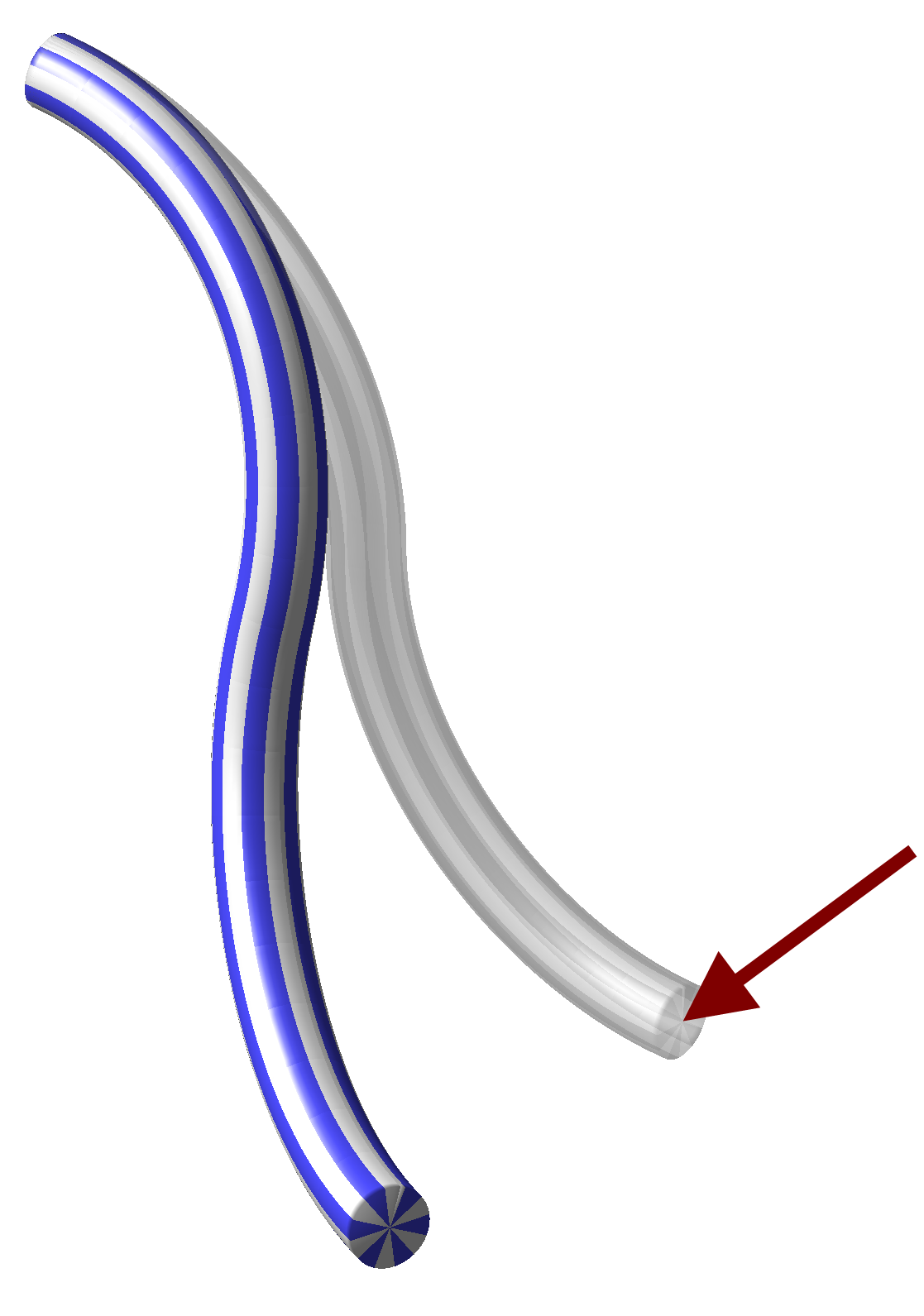}
\label{fig:sshape2}
}
\caption{S-shaped beam clamped at the upper end}
\end{figure}

\begin{figure}
\centering
\subfigure[Point load directed to the left]{
\includegraphics[width=0.70\linewidth]{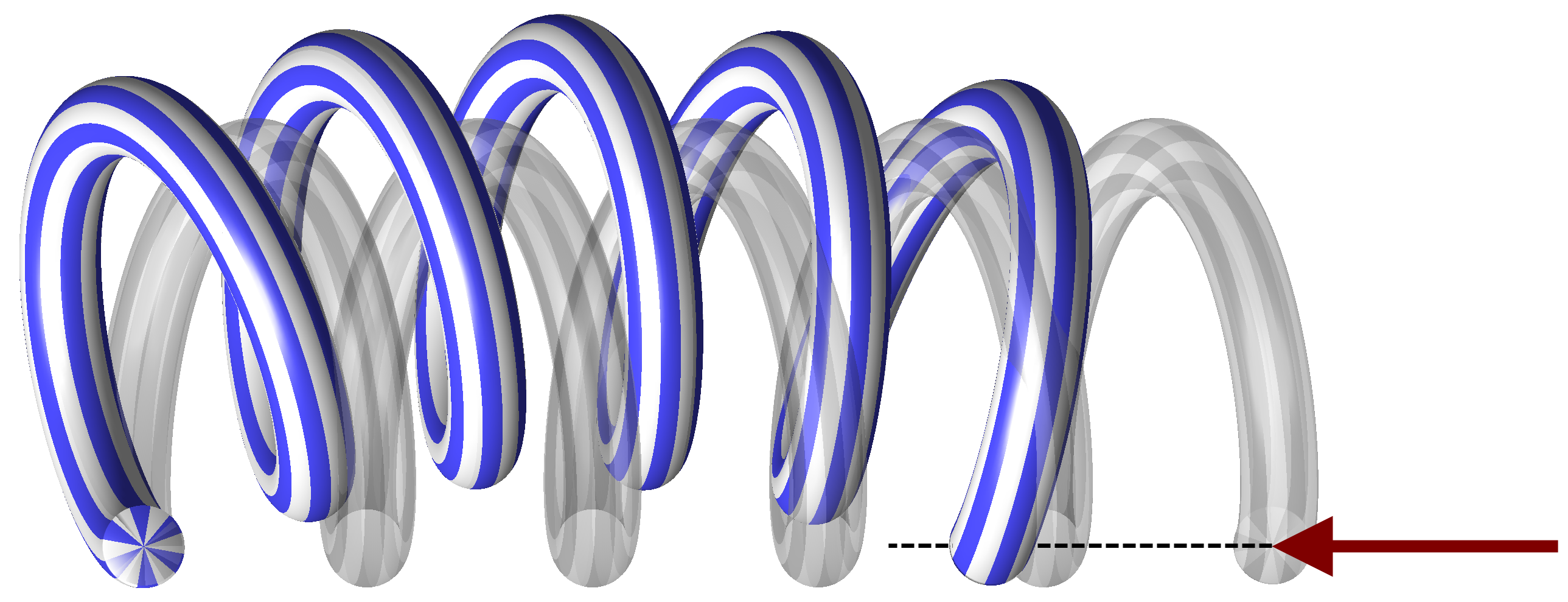}
\label{fig:spring1}
}
\subfigure[Point load directed to the right]{
\includegraphics[width=0.80\linewidth]{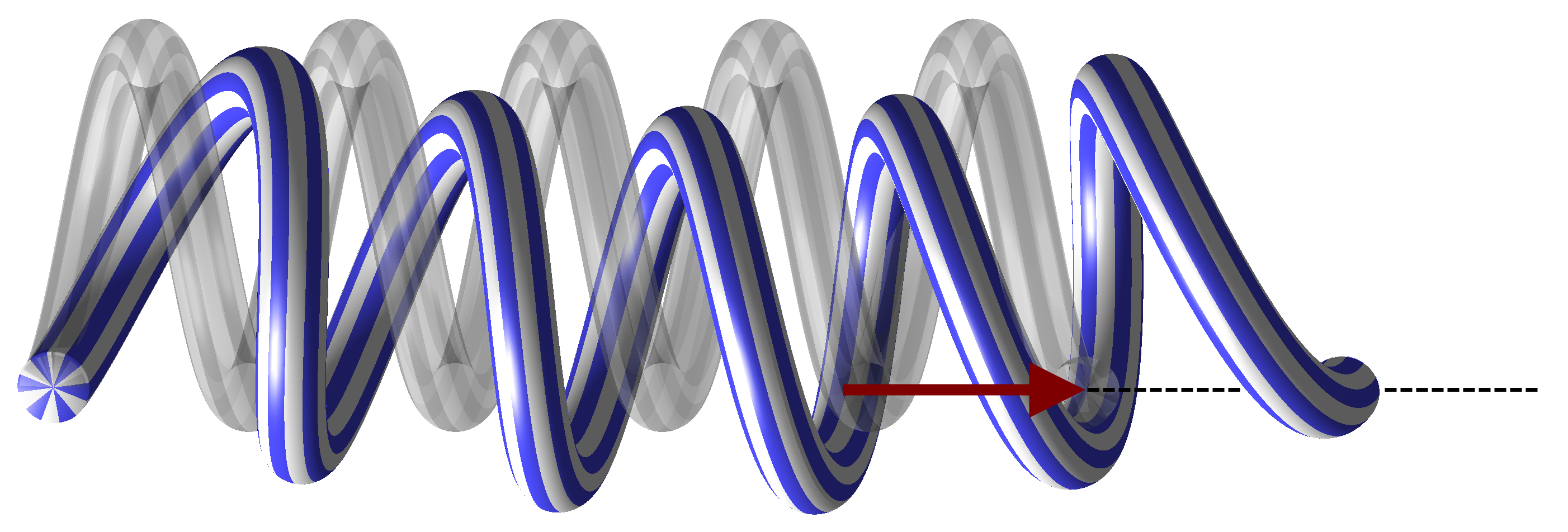}
\label{fig:spring2}
}
\caption{Helix shaped beam under point loads with the left end of the beam pinned and the right end constrained to move only in the indicated direction}
\label{fig:spring}
\end{figure}

To make it easier to explain coupling effects stemming from the initial curvature of $\Sigma$ and
also easier to compare with other beam formulations we below supply alternative expressions for the 
forces and moments \eqref{normalforce}-\eqref{torsion} where we have separated the tangential and normal plane components
of $\bfu_{\text{mid}}$ and $\bftheta$. These expressions are
\begin{align}
\bfN_\Sigma 
&=
E |A| \left( \bftcdotnabla u_t - (\bfQ_\Sigma \bfu_{\text{mid}}) \cdot \bfkappa \right) \bft
\\
\bfS_\Sigma 
&=
G |A| ( \bfQ_\Sigma \bftcdotnabla(\bfQ_\Sigma \bfu_{\text{mid}}) - (\bfQ_\Sigma \bftheta) \times\bft + u_t \bfkappa )
\label{Scurve}
\\
\bfM_\Sigma 
&=
E \bfI_\Sigma \left( \bftcdotnabla (\bfQ_\Sigma \bftheta) + \theta_t \bfkappa \right)
\label{Mcurve}
\\
\bfT_\Sigma 
&= 
G J_\Sigma \left( \bftcdotnabla\theta_t - (\bfQ_\Sigma \bftheta)\cdot\bfkappa \right) \bft
\label{Tcurve}
\end{align}
and we give some notes on how they are derived in Appendix \ref{forcesandmomentsseparation}.
The last term in each expression involves the curvature vector $\bfkappa$ and explains the
coupling effects stemming from the initial curvature of the beam. For an initially straight
beam these terms vanish and the above reduce to the expressions of classical beam theory \cite{ZiTa05}.

To give some illustration to the effects of the curvature coupling we here give a few numerical examples.
As in previous figures the beam in the undeformed state is rendered in gray and
the deformed beam is rendered in blue and white.

Applying a torque to the end of a s-shaped beam as in Figure~\ref{fig:sshape} we by \eqref{Tcurve} get a coupling
between how the twist change, i.e. $\bftcdotnabla \theta_t$, and the bending angle $\bfQ_\Sigma \bftheta$
in the direction parallel to the curvature vector $\bfkappa$.
The bending angle in turn couples to the change in midline deflection by \eqref{Scurve}, where
the change is in a direction orthogonal to the tangent. This last coupling is however
not due to curvature. Further, the coupling in \eqref{Mcurve} is also involved.

In Figure~\ref{fig:sshape2} we see the limitations of the linear elastic model.
Since the load is orthogonal to the curvature along the beam there are
no curvature coupling effects from \eqref{Scurve} coming into play. As the load is also orthogonal to
the tangent this actually becomes a pure bending problem.

While the expressions above are quite intricate as the curvature couples all the four equations,
the resulting effects seem to
be plausible as in Figure~\ref{fig:spring} where we push and pull a helix shaped spring
or as the previously mentioned example in Figure~\ref{fig:sshape}.

\section{Conclusions} \label{section:conclusion}
To our knowledge, with the exception of continuum based models \cite{Belytschko2000}, this is the first variational formulation for curved beams which is natively derived and expressed in global Cartesian coordinates. While it clearly is intuitive to derive formulations for curved beams in a local coordinate system, such as the Frenet frame, an advantage in our approach is that there is no explicit need for well defined normal directions along the beam.

There are several additional advantages of working in global coordinates.
Implementation of finite elements becomes very straightforward as it
easy to work with polynomial interpolation with global degrees of freedom.
It also gives rise to new possibilities for analysis of beam elements.
Further, global degrees of freedom are advantageous when connecting various structures.

\bibliographystyle{abbrv}
\bibliography{LapBel}

\clearpage
\appendix

\section{Expressions for strains, forces, and moments} \label{app:strain}

Note that we in some calculations below use $\bfQ_\Sigma = \bfn_1 \otimes \bfn_1 + \bfn_2 \otimes \bfn_2$
where $\bfn_1$ and $\bfn_2$ are orthogonal normals to $\Sigma$. As these are arbitrary and just used
to simplify the arguments, their use does not imply a need to define these normal directions along $\Sigma$
for the expressions below to hold.

\subsection{Strain expressions}
In this section we give derivations of the strain expressions presented
in Section~\ref{section:strain}.
Inserting the assumed displacement field in the in-line strain we have
\begin{align}
\bfeps_\Sigma^P (\bfu ) \cdot \bft = \bfeps_\Sigma^P ( \bfu_{\text{mid}} ) \cdot\bft + \bfeps_\Sigma^P (\bftheta \times \bfzeta ) \cdot \bft
\end{align}
Evaluating these expressions gives terms
\begin{align}
\bfeps_\Sigma^P ( \bfu_{\text{mid}} ) \cdot\bft
&=
\bfP_\Sigma \left(\nabla \bfu_{\text{mid}}\right)^\text{T} \cdot \bft
\\&=
\bfP_\Sigma \bftcdotnabla \bfu_{\text{mid}}
\end{align}
and
\begin{align}
\bfeps_\Sigma^P ( \bftheta \times \bfzeta ) \cdot \bft
&=
\bfP_\Sigma \left(\nabla \left(\bftheta \times \bfzeta\right)\right)^\text{T} \cdot \bft
\\&=
\bfP_\Sigma \bftcdotnabla\left(\bftheta \times \bfzeta\right)
\\&
\begin{aligned}
&=
\bfP_\Sigma \left(\bftcdotnabla\bftheta\right) \times \bfzeta
\\&\quad+
\bfP_\Sigma \bftheta \times \left(\bftcdotnabla\bfzeta\right)
\end{aligned}
\\&=
\bfP_\Sigma \left(\bftcdotnabla\bftheta\right) \times \bfzeta
\end{align}
where we use $\bftcdotnabla\bfzeta=\bfzero$.

Inserting the assumed displacement field in the shear strain we have
\begin{align}
\bfeps_\Sigma^S (\bfu ) \cdot \bft = \bfeps_\Sigma^S ( \bfu_{\text{mid}} ) \cdot\bft + \bfeps_\Sigma^S (\bftheta \times \bfzeta ) \cdot \bft
\label{someappendixeq123}
\end{align}
Evaluating the first expression gives terms
\begin{align}
\bfQ_\Sigma \left(\nabla \bfu_{\text{mid}}\right)^\text{T} \cdot \bft
&= \bfQ_\Sigma \bftcdotnabla \bfu_{\text{mid}}
\\
\bfQ_\Sigma \nabla \bfu_{\text{mid}} \cdot \bft 
&= \bfzero
\end{align}
where we use that $\bfu_{\text{mid}}$ does not change in the normal plane.
Thus, we can write
\begin{align}
2 \bfeps_\Sigma^S ( \bfu_{\text{mid}} ) \cdot \bft
&=
\bfQ_\Sigma \bftcdotnabla \bfu_{\text{mid}}
\end{align}
Now turning to evaluating the second expression of \eqref{someappendixeq123} we have the term
\begin{align}
\bfQ_\Sigma \nabla \left(\bftheta \times \bfzeta \right)^\text{T} \cdot \bft 
&=
\bfQ_\Sigma \bftcdotnabla \left(\bftheta \times \bfzeta \right) 
\\ &
\begin{aligned}
&=
\bfQ_\Sigma \left( \bftcdotnabla\bftheta \right) \times \bfzeta
\\&\quad+ \bfQ_\Sigma \bftheta \times \left( \bftcdotnabla \bfzeta \right)
\end{aligned}
\\ &=
\bfQ_\Sigma \left( \bftcdotnabla\bftheta \right) \times \bfzeta
\end{align}
where we use $\bftcdotnabla \bfzeta = \bfzero$ and we also have the term
\begin{align}
\bfQ_\Sigma \nabla (\bftheta \times \bfzeta ) \cdot \bft
&=
\sum_i \bfn_i \otimes \bfn_i \left( \nabla \left(\bftheta \times \bfzeta \right) \cdot \bft \right)
\\&=
\sum_i \bfn_i \left( \left( \bfn_i\cdot\nabla \right) \left(\bftheta \times \bfzeta \right) \cdot \bft \right)
\\&
\begin{aligned}
=
\sum_i \bfn_i &(
\left(\left( \bfn_i\cdot\nabla \right) \bftheta \times \bfzeta \right) \cdot \bft
\\&+ \left(\bftheta \times \left( \bfn_i\cdot\nabla \right) \bfzeta \right) \cdot \bft
)
\end{aligned}
\\&=
\sum_i \bfn_i \left(
\bftheta \times \bfn_i \cdot \bft
\right)
\\&=
-\sum_i \bfn_i \cdot \left(
\bfn_i \bftheta \times \bft
\right)
\\&=
-\bfQ_\Sigma \left( \bftheta\times\bft \right)
\\&=
-\bftheta\times\bft
\end{align}
where we use that $\bftheta$ does not change in the normal direction and \eqref{ndell}, i.e. $\left( \bfn_i\cdot\nabla \right) \bfzeta = \bfn_i$.
Thus, we can write
\begin{align}
2 \bfeps_\Sigma^S (\bftheta \times \bfzeta ) \cdot \bft
&=
\left( \bfP_\Sigma \bftcdotnabla \bftheta \right) \times \bfzeta
-\bftheta\times\bft
\end{align}

\subsection{Forces and moments} \label{forcesandmomentsderivation}
In this section we supply derivations of the expressions for forces and moments presented in Section~\ref{forcesandmoments}.
The definitions for in-line normal force $\bfN_\Sigma$, shear force $\bfS_\Sigma$, bending moment $\bfM_\Sigma$, and torsion $\bfT_\Sigma$ are
\begin{align}
\bfN_\Sigma
&:=
\bfP_\Sigma \int_{A} \bfsig(\bfu)\cdot\bft\, dA
\\
\bfS_\Sigma
&:=
\bfQ_\Sigma \int_{A} \bfsig(\bfu)\cdot\bft\, dA
\\
\bfM_\Sigma
&:=
\bfQ_\Sigma \int_{A}\bfzeta\times \left(\bfsig(\bfu)\cdot\bft\right) \, dA
\\
\bfT_\Sigma
&:=
\bfP_\Sigma \int_{A}\bfzeta\times \left(\bfsig(\bfu)\cdot\bft\right) \, dA
\end{align}

By \eqref{constitutive} we that
$\bfP_\Sigma \bfsig(\bfu)\cdot\bft = E \bfeps^P_\Sigma(\bfu)\cdot\bft$ and
$\bfQ_\Sigma \bfsig(\bfu)\cdot\bft = 2 G \bfeps^S_\Sigma(\bfu)\cdot\bft$.
Further, as $\Sigma$ goes through
the center of mass of $A$ we have
$\int_A \bfzeta \, dA = \bfzero$ which will be used when we evaluate the
expressions below.
Also recall \eqref{area}-\eqref{polar}, i.e. the definitions for the area, the tensor of area moments of inertia, and the polar inertia.

Evaluating the in-line normal force by using \eqref{ePu} and \eqref{ePth} we have
\begin{align}
\bfN_\Sigma
&=
\bfP_\Sigma \int_{A} \bfsig(\bfu)\cdot\bft \, dA
=
\int_{A} \bfP_\Sigma \bfsig(\bfu)\cdot\bft \, dA
\\&=
\int_{A} E \bfeps^P_\Sigma(\bfu)\cdot\bft \, dA
\\&
\begin{aligned}
=
E \int_{A} &\bfP_\Sigma \bftcdotnabla \bfu_{\text{mid}}
\\&+
\bfP_\Sigma \left( \bftcdotnabla \bftheta \times \bfzeta \right) \, dA
\end{aligned}
\label{someterm5436}
\\&=
E |A| \bfP_\Sigma \bftcdotnabla \bfu_{\text{mid}}
\end{align}
where the second term in \eqref{someterm5436} vanishes 
as $\int_{A} \bfzeta \, dA = \bfzero$.
For the shear force we use \eqref{eSu} and \eqref{eSth} which gives
\begin{align}
\bfS_\Sigma
&=
\bfQ_\Sigma \int_{A} \bfsig(\bfu)\cdot\bft \, dA
=
\int_{A} \bfQ_\Sigma \bfsig(\bfu)\cdot\bft \, dA
\\&=
\int_{A} 2 G \bfeps^S_\Sigma(\bfu)\cdot\bft \, dA
\\ &
\begin{aligned}
=
G \int_{A}
&\bfQ_\Sigma \bftcdotnabla \bfu_{\text{mid}}
\\&+ \bfQ_\Sigma \left( \bftcdotnabla \bftheta \times \bfzeta \right)
- \bftheta\times\bft
\, dA
\end{aligned}
\label{someterm9483}
\\ &=
G |A| \left( \bfQ_\Sigma \bftcdotnabla \bfu_{\text{mid}} - \bftheta\times\bft \right)
\end{align}
where the middle term in \eqref{someterm9483} vanishes $\int_A \bfzeta \, dA = \bfzero$.
We evaluate the bending moment by using \eqref{ePu} and \eqref{ePth} which gives
\begin{align}
\bfM_\Sigma
&=
\bfQ_\Sigma \int_{A}\bfzeta\times \left(\bfsig(\bfu)\cdot\bft\right) \, dA
\\&=
\int_{A}\bfzeta \times \left(\bfP_\Sigma \bfsig(\bfu)\cdot\bft\right) \, dA
\\&=
\int_{A}\bfzeta \times \left( E \bfeps^P_\Sigma(\bfu)\cdot\bft \right) \, dA
\\&=
E \int_{A}\bfzeta \times \left( \bfP_\Sigma \bftcdotnabla \bfu_{\text{mid}} \right) \, dA
\label{someterm4956}
\\&\quad+
E \int_{A}\bfzeta \times 
\left(\left( \bfQ_\Sigma \bftcdotnabla \bftheta \right) \times \bfzeta \right)
\, dA
\\&=
E \int_{A}\bfzeta \times 
\left( \bfP_\Sigma \left( \bftcdotnabla \bftheta \times \bfzeta \right) \right)
\, dA
\\&=
E \int_{A} (\bfzeta\times\bft)\otimes(\bfzeta\times\bft) \, dA \bftcdotnabla \bftheta
\\&=
E \bfI_\Sigma
\bftcdotnabla \bftheta
\end{align}
where \eqref{someterm4956} vansishes as $\int_{A} \bfzeta \, dA = \bfzero$.
Finally, using \eqref{eSu} and \eqref{eSth} we evaluate the torsion which gives
\begin{align}
\bfT_\Sigma
&=
\bfP_\Sigma \int_{A}\bfzeta\times \left(\bfsig(\bfu)\cdot\bft\right) \, dA
\\&=
\int_{A}\bfzeta\times \left(\bfQ_\Sigma \bfsig(\bfu)\cdot\bft\right) \, dA
\\&=
G \int_{A}\bfzeta \times \left( 2\bfeps^S_\Sigma(\bfu)\cdot\bft \right) \, dA
\\&=
G \int_{A}\bfzeta \times \left( \bfQ_\Sigma \bftcdotnabla \bfu_{\text{mid}} \right) \, dA
\label{someterm4321}
\\&\quad+
G \int_{A}\bfzeta \times \left( \bfQ_\Sigma \left( \bftcdotnabla \bftheta  \times \bfzeta \right)
-\bftheta\times\bft \right) \, dA
\label{someotherterm4321}
\\&=
G \int_{A}\bfzeta\cdot\bfzeta \, dA \bfP_\Sigma \bftcdotnabla \bftheta 
\\&=
G J_\Sigma \bfP_\Sigma \bftcdotnabla \bftheta
\end{align}
where \eqref{someterm4321} and the last term in \eqref{someotherterm4321} vanishes 
as $\int_{A} \bfzeta \, dA = \bfzero$.

\subsection{Tangential separation} \label{forcesandmomentsseparation}
We may separate any vector $\bfv$ into components of the tangent and the normal plane
by $\bfv = \bfP_\Sigma \bfv + \bfQ_\Sigma \bfv = \bft v_t + \bfQ_\Sigma \bfv$.
Using the that $\bfQ_\Sigma$ is a projection operator, i.e. $\bfQ_\Sigma = \bfQ_\Sigma \bfQ_\Sigma$,
we have
\begin{align}
\bftcdotnabla \bfv
&=
\bftcdotnabla
\left( \bft v_t + \bfQ_\Sigma \bfv \right)
\\&
\begin{aligned}
&= v_t \bftcdotnabla\bft + \bft \bftcdotnabla v_t
\\&\quad+
\bftcdotnabla\left(\bfQ_\Sigma\right) \left( \bfQ_\Sigma \bfv\right)
\\&\quad+ \bfQ_\Sigma \bftcdotnabla \left(\bfQ_\Sigma \bfv \right)
\end{aligned}
\\&
\begin{aligned}
&= v_t \bfkappa + \bft \bftcdotnabla v_t
\\&\quad+
\left(-\bft\otimes\bfkappa-\bfkappa\otimes\bft \right) \bfQ_\Sigma \bfv
\\&\quad+ \bfQ_\Sigma \bftcdotnabla \left( \bfQ_\Sigma \bfv \right)
\end{aligned}
\\&
\begin{aligned}
&=
v_t \bfkappa + \bft \bftcdotnabla v_t
-\bft \left(\bfkappa\cdot\bfv\right)
\\&\quad + \bfQ_\Sigma \bftcdotnabla \left(\bfQ_\Sigma \bfv \right)
\end{aligned}
\end{align}
which in turn this gives the relationships
\begin{align}
\bfQ_\Sigma \bftcdotnabla \bfv
&=
v_t \bfkappa 
+ \bfQ_\Sigma \bftcdotnabla \left(\bfQ_\Sigma\bfv\right)
\\
\bfP_\Sigma \bftcdotnabla \bfv
&=
\bft \bftcdotnabla v_t
-\bft \left(\bfkappa\cdot\bfv\right)
\end{align}

\end{document}